\documentclass[11pt]{article}
\usepackage{color}
\usepackage{amsmath,amsthm,amssymb,amsfonts}
\usepackage{epsfig}
\usepackage{graphicx}
\usepackage{yhmath} 
\usepackage{bm}
\usepackage{dsfont}
\def\R{\mathbb R}
\def\N{\mathbb N}
\def\Z{\mathbb Z}

\def\e{\varepsilon}
\def\trait (#1) (#2) (#3){\vrule width #1pt height #2pt depth #3pt}
\def\fin{\hfill\trait (0.1) (5) (0) \trait (5) (0.1) (0) \kern-5pt
\trait (5) (5) (-4.9) \trait (0.1) (5) (0)}

\numberwithin{equation}{section}

\newcommand{\be}{\begin{equation}}
\newcommand{\ee}{\end{equation}}
\newcommand{\baa}{\begin{array}}
\newcommand{\eaa}{\end{array}}
\newcommand{\ba}{\begin{eqnarray}}
\newcommand{\ea}{\end{eqnarray}}

\newcommand{\ban}{\begin{eqnarray*}}
\newcommand{\ean}{\end{eqnarray*}}

\newcommand{\dis}{\displaystyle}
\newtheorem{theo}{\bf Theorem}[section]
\newtheorem{lem}[theo]{\bf Lemma}
\newtheorem{pro}[theo]{\bf Proposition}

\newtheorem{defi}[theo]{\bf Definition}
\newtheorem{rem}[theo]{\bf Remark}

\def\ds{\rightarrow}
\def\e{\varepsilon}

\def\Ac{{\cal A}}

\def\apg{\left\{}
\def\chg{\right\}}
\def\apt{\left(}
\def\cht{\right)}

\def\ap{\left.}
\def\ch{\right.}

\def\a{a}

\def\A{{\cal A}}

\def\eps{\varepsilon}

\def\Euno{\mathcal{E}_1}
\def\Edue{\mathcal{E}_2}
\def\Eh{\mathcal{E}_\H}

\newcommand{\mB}{\mathcal{B}}
\newcommand{\mBL}{\mathcal{BL}}
\newcommand{\cob}{\overline{\mathop{\rm co}}}
\def\ue{u_\eps}

\def\trajYy{Y_{y}}
\def\trajLy{L_{y}}

\def\uep{U^+_\eps}
\def\uem{U^-_\eps}
\def\Hbarp{\bar{H}^+}
\def\Hbarm{\bar{H}^-}
\renewcommand{\H}{\mathcal{H}}
\def\y{\frac{x}{\eps}}
\def\up{U^+}
\def\um{U^-}

\def\HT{H_T}

\def\bH{b_{\H}}

\def\lH{l_{\H}}

\def\HTreg{{ H}^{\rm reg}_T}

\newcommand{\trajxo}{X^{\e}_{x_0}}
\def\dottrajxo{\dot{X}^{\e}_{x_0}}
\def\Ta{{\cal T}^{\e}_{x_0}}
\def\Treg{{\cal T}^{\rm reg,\e}_{x_0}}

\def\upcor{{V}^{+}}
\def\umcor{{V}^{-}}
\def\Vrho{U^{\rho+}}

\def\Vrho{V^{\rho+}}

\def\Hbarp{\bar{H}^+}
\def\Hbarm{\bar{H}^-}

\def\lrm{\lambda_\rho^{-}}

\def\nor{\mathbf{n}}

\def\HHm{\mathbf{H}^-}
\def\HHp{\mathbf{H}^+}

\def\valpham{v^{\rho-}}

\def\valphap{v^{\rho+}}
\def\walphap{w^{\rho+}}

\def\ualpham{V_{\rho}^{+}}
\def\ualpham{V_{\rho}^{-}}
\def\trajyoyo{Y_{y_0}}
\def\dottrajyoyo{\dot{Y}_{y_0}}
\def\Tayo{{\cal T}_{y_0}}
\def\Tayoreg{{\cal T}^{\rm reg}_{y_0}}

\def\Taxbarrareg{{\cal T}^{\rm reg}_{\bar{x}}}

\def\xb{{\bar{x}}}

\def\Tayarrareg{{\cal T}^{\rm reg}_{y}}
\definecolor{NT}{gray}{0.60}

\begin{document}

\title{\bf   Homogenization Results for a Deterministic Multi-domains Periodic Control Problem}
\author{G.Barles, A. Briani, E. Chasseigne, 
\thanks{Laboratoire de Math\'ematiques et Physique Th\'eorique (UMR CNRS 7350), F\'ed\'eration Denis Poisson (FR CNRS 2964),
  Universit\'e Fran\c{c}ois Rabelais, Parc de Grandmont,
 37200 Tours, France. email: Guy.Barles@lmpt.univ-tours.fr,  ariela.briani@lmpt.univ-tours.fr, Emmanuel.Chasseigne@lmpt.univ-tours.fr}     
 N. Tchou,   \thanks{ IRMAR, UMR 6625,  UniversitŽ de Rennes 1, Rennes, France. email: nicoletta.tchou@univ-rennes1.fr}
  }
\maketitle

\begin{abstract}
We consider homogenization problems in the framework of deterministic
optimal control when the dynamics and running costs are completely
different in two (or more) complementary domains of the space $\R^N$.
For such optimal control problems, the three first authors have shown
that several value functions can be defined, depending, in particular,
of the choice is to use only ``regular strategies'' or to use also
``singular strategies''. We study the homogenization problem in these two
different cases. It is worth pointing out that, if the second one can be
handled by usual partial differential equations method " \`a la
Lions-Papanicolaou-Varadhan" with suitable adaptations, the first case
has to be treated by control methods (dynamic programming).
\end{abstract}

\\
{\bf AMS Class. No}:
49L20,   
49L25.   
35F21   
93C70   
35B25   

\section{Introduction}

In order to describe the homogenization problems we address in this article, we consider a partition of $\R^N$ as $\Omega_1\cup\Omega_2\cup\H$ where $\Omega_1, \Omega_2$ are open subsets of $\R^N$, $\Omega_1\cap\Omega_2=\emptyset$ and $\H=\partial\Omega_1=\partial\Omega_2$. We assume that the $\Omega_i$'s are $\Z^N$-periodic, i.e. $x+z
\in \Omega_i$ for all $x\in \Omega_i$ and $z\in \Z^N$.

The homogenization problems can be written from the partial differential equations (pde in short) point of view as
\begin{equation}\label{Bellman-Om}
\begin{array}{cc}
{\dis \lambda \ue(x) + H_1(x,\y,D\ue(x)) = 0     }&   \hbox{ in   }  \e   \Omega_1\,, \\
& \\
{\dis \lambda \ue(x)  + H_2(x,\y,D\ue(x)) = 0     }&     \hbox{ in   } \e   \Omega_2\,,
\end{array}
\end{equation}
where $\e$ is a small positive parameter which is devoted to tend to 0, the actualization factor $\lambda$ is positive and $H_1, H_2$ are classical Hamiltonians of deterministic control problems, which are of the form ($i=1,2$)
\be  \label{def:Ham}
H_i(x,y,p):=\sup_{ \alpha_i \in A_i} \apg  -b_i(x,y,\alpha_i) \cdot p - l_i(x,y, \alpha_i) \chg\,,   \quad \hbox{for  }x \in \R^N  \:, y \in \overline{\Omega}_i  \:, p \in \R^N\, .
\ee
Precise assumptions will be given later on but we already mention that the functions $b_i$ and $l_i$ satisfy the most classical regularity and boundedness assumptions and $b_i(x, y, \alpha_i)$ and  $l_i(x , y, \alpha_i)$ are $ \Z^N
$-periodic in $y$, for any $x$ and $\alpha_i$.

Of course, Equations~\eqref{Bellman-Om} have to be completed by suitable conditions on the hypersurface $\H$\footnote{Our assumptions below actually imply that $\H$ is a $W^{2,\infty}$-hypersurface} and this was the purpose of \cite{BBC,BBC2} to see what kind of conditions have to be imposed. Unfortunately the classical Ishii inequalities 
\begin{equation}\label{Bellman-H-sub}
\min\{   \lambda \ue + H_1(x,\y,D\ue)    ,  \lambda \ue  + H_2(x,\y,D\ue)  \}\leq 0  \quad\hbox{on   } \e   \H \; ,
\end{equation}
and 
\begin{equation}\label{Bellman-H-sup}
 \max\{ \lambda \ue + H_1(x,\y,D\ue)  ,     \lambda \ue  + H_2(x,\y,D\ue)  \}\geq 0   \quad\hbox{on   }  \e  \H \; ,
\end{equation}
are not sufficient to have a well-posed problem and \eqref{Bellman-Om}-\eqref{Bellman-H-sub}-\eqref{Bellman-H-sup} has a maximal solution denoted by $\uep$ and a minimal solution denoted by $\uem$ which can both be described in terms of control. We refer the reader to Section~\ref{sec:ueps} for a complete description of the control problems for $\uep$ and $\uem$ but we just mention that, while $\uep$ is built by using so-called ``regular strategies'', $\uem$ is built by using all kind of strategies and in particular ``singular strategies'' which are excluded in the case of $\uep$.

We want to describe the asymptotic behavior as $\e\rightarrow 0$ of the maximal solution $\uep$ and the minimal solution $\uem$.
The results in \cite{BBC,BBC2} imply that $\uem$ can be characterized through pdes, by adding a suitable subsolution condition on $\H$, while this is not the case anymore for $\uep$ which is just the maximal subsolution of \eqref{Bellman-Om}-\eqref{Bellman-H-sub}-\eqref{Bellman-H-sup}. 
The consequence for our study is immediate: while for $\uem$ we can
follow and adapt the classical pde arguments of Lions, Papanicolaou \&
Varadhan \cite{LPV} and Alvarez \& Bardi \cite{AB1,AB2}, this is not the
case anymore for $\uep$ where even if we follow closely the pde ideas,
we have to perform all the argument on the control formulas. In that way
we are close to some of the arguments of the weak KAM theory (see Fathi \cite{Fa1, Fa2, Fa3}).

For the convergence of $\uep$, some specific technical difficulties appear which are solved by an approximation of the cell problem : such ideas, in a slightly different context, are already used in Barles, Da Lio, Lions \& Souganidis \cite{BDLLS} (see also \cite{B1}).

We point out that Forcadel and Rao \cite{FR} studied such homogenization problems in a multi-domain framework : they are able to treat cases where the boundaries of the $\Omega_i$ are not smooth but for problems set in $\R^2$ and only in the $\uem$ case.  
As related works we mention  the study of Hamilton-Jacobi equations on networks  \cite{ACCT,IM,IMZ}. 

The article is organized as follows : in Section~\ref{solvisco}, we recall basic facts and stability results for the pde approach of \eqref{Bellman-Om}-\eqref{Bellman-H-sub}-\eqref{Bellman-H-sup}. In Section~\ref{sec:ueps}, we describe the control problems for $\uep$ and $\uem$ : we give precise definitions of ``regular strategies'' and ``singular strategies''. The next two sections are devoted to the study of the homogenization problems for $\uem$ and $\uep$ respectively: we follow a (rather) classical double-scale approach by first studying the cell problems and then we use the solutions of the cell problems to deduce the convergence. 
In Section ~\ref{sect:oneD} we give an explicit example in dimension 1 where the effective Hamiltonians describing the asymptotic behavior  of $\uem$ and $\uep$, are different.
Finally in the appendix we provide several technical results which are useful for the convergence proofs.

\section{ Different notions of viscosity solutions for multi-domains problems}  \label{solvisco}

This section is devoted to  the description of the precise  definition of  viscosity solutions for problems like \eqref{Bellman-Om}-\eqref{Bellman-H-sub}-\eqref{Bellman-H-sup}.  
For the introduction and all the details on these definitions we refer to \cite{BBC,BBC2}  and the reference therein. 

Let us remember that in this paper we are considering a partition of $\R^N$ as $\Omega_1\cup\Omega_2\cup\H$ where $\Omega_1, \Omega_2$ are open subsets of $\R^N$, $\Omega_1\cap\Omega_2=\emptyset$ and $\H=\partial\Omega_1=\partial\Omega_2$ is a regular hypersurface ($W^{2,\infty}$).
For $y\in\H$ we denote by $T_y\H$, the tangent space to $\H$ at $y$
and $\dis <\cdot,\cdot>_{T_y\H}$ is the scalar product in this tangent space. 

We consider the   general  function  $\mathbf{G} : \R^N \times \R^N \times \R^N \mapsto \R$  that can be differently defined on $\Omega_1$, $\Omega_2$ and $\H$ by
$$
\mathbf{G}(x,y,p ): =\apg \begin{array}{cl}
{\dis  G_1(x,y,p)    }&   \hbox{ if   } y \in  \Omega_1\,, \\
{\dis  G_\H(x,y,p)    }&   \hbox{ if   } y \in  \H  \hbox{ and   } p\in T_y\H\,, \\
{\dis  G_2(x,y,p)    }&     \hbox{ if  }   y \in  \Omega_2\,.
\end{array} \ch 
$$
We define also 
$$
G(x,y,p): =\apg \begin{array}{cc}
{\dis  G_1(x,y,p)    }&   \hbox{ if   } y \in  \Omega_1\,, \\
{\dis  G_2(x,y,p)    }&     \hbox{ if  }   y \in  \Omega_2\,.
\end{array} \ch 
$$
Note that to be consistent with the homogenization problem  we will always assume that  each $G_i$ is not only  defined in $\Omega_i$ but in all $\R^N$.   \\
First of all we recall the classical H. Ishii definition of  discontinuous viscosity solution for a discontinuous Hamiltonian $G$ (see \cite{HI} and also \cite{BaCIME}).
Given a real number  $\rho \geq 0$ and a function $f : \R^N \times \R^N\mapsto \R$ a viscosity solution of problem
\be \label{progenNOH}
\rho u(y)+ G(x,y, Du(y)) = f(x,y)  \hbox{ in }  \R^N
\ee
is defined as follows.
\begin{defi} \label{solviscogenNOH}  
We say that a bounded usc function $u$ is a 
\underline{{\rm subsolution} of   \eqref{progenNOH}} if it verifies the following inequalities in the viscosity sense 
\be
\label{Ishiisub}
  \rho u(y) + G_1(x,y,Du(y)) \leq  f(x,y)   \;   \mbox{ for }  y  \in \Omega_1    \:,   \;    \rho u(y) + G_2(x,y,Du(y)) \leq  f(x,y)      \mbox{ for }  y  \in \Omega_2   \:, 
\ee
\be
\label{IshiisubH}
   \rho u(y) + \min\{G_1(x,y,Du(y)), G_2(x,y,Du(y)) \} \leq  f(x,y)     \mbox{ for }  y  \in \H .
\ee

We say that a  lsc function $v$ is a  \underline{{\rm supersolution} of   \eqref{progenNOH}} 
if it verifies the following inequalities in the viscosity sense 
\be
\label{Ishiisup}
 \;   \rho u(y) + G_1(x,y,Du(y)) \geq  f(x,y)     \mbox{ for }  y  \in \Omega_1  \:,   \quad    
 \quad   \;  \rho u(y) + G_2(x,y,Du(y)) \geq  f(x,y)    \mbox{ for }  y  \in \Omega_2  \:,  
\ee
\be
\label{IshiisupH}
     \rho u(y) + \max\{G_1(x,y,Du(y)), G_2(x,y,Du(y)) \} \geq  f(x,y)  \mbox{ for }  y  \in \H  .
\ee

\
We say that a  bounded continuous  function $w$ is a  \underline{{\rm solution} of   \eqref{progenNOH}}  it is both a sub and a supersolution. 
\end{defi}

We give now a new definition of solution taking in account  a tangential equation on $\H$ of the following problem:
\be \label{progen}
\rho u(y)+ \mathbf{G}(x,y, Du(y)) = f(x,y)  \hbox{ in }  \R^N.
\ee
The main difference is the following definition of viscosity subsolution  for the tangential Hamiltonian on $\H$:
\begin{defi} \label{defi:sousolH}
 An  usc  function $u : \H\to\R $ is a viscosity subsolution of 
  $$ \rho u(y) + G_\H(x,y,D_\H u(y))=f(x,y)\quad \text{on   } \H$$   
     if, for any $\phi\in C^1(\R^N)$ and any maximum point
    $y$ of $z\mapsto u(z)-\phi(z)$ in $\H $,
    one has
    $$ \rho u(y)+ G_\H (x,y,D_\H\phi(y) )\leq f(x,y)\; ,$$ 
    where $D_\H\phi(y)$ means 
    the gradient of the restriction of $\phi$  to $\H$ at  point $y$
    (which belongs to $T_y\H$)\footnote{Note that, if $\nor (y)$ is a
	unitary normal vector to $\H$ at $y$, then
 $D_{\H} \phi (y)=D\phi(y) -(D\phi(y) \cdot \nor(y)) \nor(y)$.}.
\end{defi}
We define the viscosity solutions for multi-domain problems as Ishii classical solutions adding for subsolutions the previous condition, more precisely :
\begin{defi}   \label{solviscogen}  
We say that a bounded usc function $u$ is a 
\underline{{\rm subsolution} of   \eqref{progen}} if it verifies  \eqref{Ishiisub},  \eqref{IshiisubH} and
\be    \rho u(y) +  G_\H(x,y, D_\H u(y)) \leq  f(x,y)   \mbox{ in the sense of Definition  \ref{defi:sousolH} above},  \mbox{ for }  y  \in \H.
\ee
We say that a lsc function $v$ is a  \underline{{\rm supersolution} of   \eqref{progen}} 
if it verifies  \eqref{Ishiisup} and  \eqref{IshiisupH}.



\
We say that a  bounded continuous  function $w$ is a  \underline{{\rm solution} of   \eqref{progen}}  it is both a sub and a supersolution. 
\end{defi}

Note that an analogous definition can be given for a solution $u$ depending on the $x$-variable instead of  the $y$-variable, so we will not detail it.


Moreover let us remark that the function $G_\H (x,y, p)$ doen't need to be defined for any $p\in\R^N$ but only in the tangent space $T_y\H$.



Now we focus our attention on the Bellman Equations when we have Hamiltonians of the form \eqref{def:Ham}. We state our main assumptions and define the different type of dynamics and Hamiltonians on the interfaces. 

\begin{itemize}

\item[{[H0]}]  For $i=1,2$, $A_i$ is a compact metric space and $b_i : \R^N \times \R^N \times A_i  \ds \R^N$ is a continuous bounded function. 
More precisely, there exists  
$M_b >0$, such that for any $x \in \R^N$, $y \in \R^N$ and $\alpha_i \in A_i$, $i=1,2$,
$$ |b_i(x,y,\alpha_i) | \leq M_b \; . $$ 
For any $x \in \R^N$ and  $\alpha_i \in A_i$, the function $b_i(x,\cdot, \alpha_i)$ is $\Z^N$-periodic. \\
 Moreover, there exists $L_i \in \R$ such that, for any $x,z \in \R^N$, $y \in \R$ and $\alpha_i \in A_i$
$$ |b_i(x,y,\alpha_i)-b_i(z,y,\alpha_i)|\leq L_i |x-z|\; , $$ 
and there exists $\bar{L}_i \in \R$ such that, for any $x\in \R^N$, $y,w \in \R$ and $\alpha_i \in A_i$
$$ |b_i(x,y,\alpha_i)-b_i(x,w,\alpha_i)|\leq \bar{L}_i |y-w|\;.$$ 

\item[{[H1]}]  For $i=1,2$, the function $l_i : \R^N \times \R^N \times A_i  \ds \R^N$ is a continuous, bounded function.  More precisely, there exists  
$M_l >0$, such that for any $x \in \R^N$, $y \in \R^N$ and $\alpha_i \in A_i$, $i=1,2$,
$$ |l_i(x,y,\alpha_i) | \leq M_l \; . $$ 
Moreover,   for any $x \in \R^N$, $\alpha_i \in A_i$,  the function $l_i(x,\cdot, \alpha_i)$ is $(2\Z)^N$-periodic. \\
There exists a modulus $\omega_l(\cdot)$ such that, for any $x,z \in \R^N$, $y \in \R^N$ and $\alpha_i \in A_i$
$$ |l_i(x,y,\alpha_i)-l_i(z,y,\alpha_i)|  \leq \omega_{l}( x-z)\;\; , $$ 
and there exists $\bar{L}_{i,l}\in \R$ such that, for any $x\in \R^N$, $y,w \in \R^N$ and $\alpha_i \in A_i$
$$ |l_{i}(x,y,\alpha_i)-l_{i}(x,w,\alpha_i)|\leq \bar{L}_{i,l} |y-w|\;.$$ 


%
\end{itemize}

\noindent Let us recall that a modulus is a function $\psi:\R^N\rightarrow\R^+$ such that $\psi(z)=\omega(\vert z \vert)$ where
$\omega$ is an increasing function $\omega:\R^+\rightarrow \R^+$ 
such that $\displaystyle\lim_{t\rightarrow 0+}\omega(t)=0$.

\begin{itemize}
\item[{[H2]}] For each $x,y  \in \R^N$ and $i=1,2$ the sets $\displaystyle \cup_{\alpha_i \in A_i } (b_i(x,y, \alpha_i),l_i(x, y,\alpha_i))  $,  are closed and convex. There is a $\delta>0$ such that for any
	$i=1,2$, $x  \in \R^\N$, \:  $y \in \H$. 
    \begin{equation}\label{cont-ass}
        \mathbf{B}_i(x,y)\supset \{|z|\leq\delta\}
     \end{equation} 
     where $\mathbf{B}_i(x,y):= \big\{ b_i(x,y,\alpha_i) : \alpha_i \in
	 A_i \big\}$\,.  
\end{itemize}
\noindent We set 
$$A:=\{(\alpha_1, \alpha_2, \mu);\alpha_i \in A_i,\mu \in [0,1]\}.$$
For $x\in\R^N$, $y \in \H$, $a=(\alpha_1, \alpha_2, \mu)\in A$ we denote by
  \begin{equation}
  \label{bH}
  \bH (x,y,a):=\mu b_1 (x,y,\alpha_1) + (1-\mu)b_2(x,y,\alpha_2)\,,
   \end{equation} 
and
  \begin{equation}\label{lH}
  \lH (x,y,a):=\mu l_1 (x,y,\alpha_1) + (1-\mu)l_2(x,y,\alpha_2)\,. 
  \end{equation}

\noindent The set of tangential controls is given by:
 $$A_0(x,y):=\{ a\in A : \bH(x,y,a)\cdot  \nor_1(y) = 0\}$$
 where $ \nor_i(y)$ is the unitary normal exterior vector to $\Omega_i$ in $y$, and the subset of $A_0(x,y)$ of "regular" tangential controls is given by
$$A^{\rm reg}_0(x,y):=\{ a\in A_0(x,y) : b_i (x,y,\alpha_i)\cdot  \nor_i(y) \geq 0\},$$ 
the tangential Hamiltonians:
\begin{equation}  \label{def:HamHT}
\HT(x,y,p):=\sup_{a\in A_0(x,y)}
 \big\{  - <\bH(x,y,\a),p>_{T_y\H}   -  \lH(x,y,\a)   \big\},
\end{equation}
\begin{equation}  \label{def:HamHTreg}
\HTreg(x,y,p):=\sup_{a\in A^{reg}_0(x,y)} \big\{   - <\bH(x,y,\a),p>_{T_y\H}  -  \lH(x,y,\a)   \big\},
\end{equation}
where $p\in T_y\H$ and  $\bH(x,y,\a)$ has been identified with its orthogonal projection on $T_y\H$ 
and the Hamiltonians associated with $H_1$ on $\Omega_1$, $H_2$ on $\Omega_2$ (defined in  (\ref{def:Ham})) and respectively $\HT$ and $\HTreg$ on $\H$:
$$
\label{def:Hamgras}
\mathbf{H^-}(x,y,p ): =\apg \begin{array}{cl}
{\dis  H_1(x,y,p)    }&   \hbox{ if   } y \in  \Omega_1\,, \\
{\dis  \HT(x,y,p)    }&   \hbox{ if   } y \in  \H  \hbox{ and   } p\in T_y\H\,, \\
{\dis  H_2(x,y,p)    }&     \hbox{ if  }   y \in  \Omega_2\,.
\end{array} \ch 
$$
$$
\label{def:Hamgrasreg}
\mathbf{H}^{^+}(x,y,p ): =\apg \begin{array}{cl}
{\dis  H_1(x,y,p)    }&   \hbox{ if   } y \in  \Omega_1\,, \\
{\dis  \HTreg(x,y,p)    }&   \hbox{ if   } y \in  \H  \hbox{ and   } p\in T_y\H\,, \\
{\dis  H_2(x,y,p)    }&     \hbox{ if  }   y \in  \Omega_2\,.
\end{array} \ch 
$$
$$
\label{def:Hamnorm}
H(x,y,p): =\apg \begin{array}{cc}
{\dis  H_1(x,y,p)    }&   \hbox{ if   } y \in  \Omega_1\,, \\
{\dis  H_2(x,y,p)    }&     \hbox{ if  }   y \in  \Omega_2\,.
\end{array} \ch 
$$

\noindent We are interested in the following equations:

\begin{equation}\label{Bellman-OmH+}
\dis \lambda v(x) + \mathbf{H^+_{\e}}(x,\y,Dv(x)) = 0  \hbox{ in }  \R^N
\end{equation}

\begin{equation}\label{Bellman-OmH-}
\dis \lambda v(x) + \mathbf{H^-_{\e}}(x,\y,Dv(x)) = 0  \hbox{ in }  \R^N
\end{equation}

\begin{equation}\label{Bellman-OmH}
\dis \lambda v(x) + H_{\e}(x,\y,Dv(x)) = 0  \hbox{ in }  \R^N
\end{equation}
where $\lambda>0$ and $ \mathbf{H^+_{\e}}(x,y,p)$ is associated with $H_1$ on $\e\Omega_1$, $H_2$ on $\e\Omega_2$ and $\HTreg$ on $\e\H$,
$ \mathbf{H^-_{\e}}(x,y,p)$ is associated with $H_1$ on $\e\Omega_1$, $H_2$ on $\e\Omega_2$ and $\HT$ on $\e\H$, $H_{\e}(x,y,p)$ is associated with $H_1$ on $\e\Omega_1$, $H_2$ on $\e\Omega_2$ (in the following we delete the index $\e$ in these notations for the sake of simplicity).

\section{Setting  the optimal control  problem at $\eps$-fixed
}\label{sec:ueps}

The aim of this section is to give  the precise definition of infinite horizon control problems  whose value functions are "solutions" of  the Hamilton-Jacobi-Bellman equations  (\ref{Bellman-OmH-}) or  (\ref{Bellman-OmH+}).  


Note that, assumptions [H0], [H1], 
are the classical hypotheses used in infinite horizon control problems. We have strengthened them in [H2] in order to have uniformly Lipschitz continuous value functions.  Let us remark also that the first part of assumption [H2] avoids the use of relaxed controls.

In order to define  the optimal control problems in all $\R^N$, we have to define the dynamics and therefore we are led to consider an ordinary differential equation with discontinuous right-hand side. This kind of ode has been treated for the first time in the pioneering work of Filippov~\cite{Fi}. We are going to define the trajectories of our optimal control problem by using the approach through differential inclusions which is rather convenient here. This  approach has been introduced in \cite{Wa} (see also \cite{AF}) and has become now classical.  To do so in a more general setting, and since the controllability condition (\ref{cont-ass}) plays no role in the definition of the dynamic, we are going to use Assumption [H2]$_{nc}$ which is [H2] without (\ref{cont-ass}).

Fix $\e >0$. Our trajectories $\trajxo(\cdot)=\big(X^\e_{x_0,1},X^\e_{x_0,2},\dots,X^\e_{x_0,N} \big)(\cdot)$ are Lipschitz continuous functions which are solutions  of the following differential inclusion
\begin{equation} \label{def:traj}
\dottrajxo (t) \in \mB\Big(\trajxo (t), \frac{\trajxo (t) }{\e}\Big)  \quad \hbox{for a.e.  } t \in (0,+\infty)  \: ; \quad \trajxo (0)=x_0
\end{equation}
where, for any $x,y \in \R^N$, 
\[
\mB(x,y):= \apg
\begin{array}{cc}
  \mathbf{B}_1(x,y)   &   \mbox{  if }  y \in \Omega_1\,   \\
  \mathbf{B}_2(x,y) &    \mbox{  if }   y \in  \Omega_2\,  \\
   \cob \big( \mathbf{B}_1(x,y) \cup \mathbf{B}_2(x,y) \big)  &        \mbox{  if }  y \in \H\
\end{array},
\ch
\]
the notation $\cob(E)$ referring to the convex closure of the set $E\subset\R^N$ and we recall that $\mathbf{B}_i(x,y)$ are defined in [H2].

We denote by $A$  the set $A:= A_1 \times A_2 \times [0,1]$ and we set $\Ac := L^\infty (0,+\infty; A)$.\\
We have the following
\begin{theo}\label{def:dynAssume} {\bf \cite[Thm. 2.1]{BBC2}}
    Assume [H0], [H1] and [H2]$_{nc}$. Fix $\e >0$, then

\noindent $(i)$ For each $x_0 \in \R^N$, there exists a Lipschitz function $\trajxo : [0,\infty[ \ds \R^N$
which is a solution of  the differential inclusion  \eqref{def:traj}.

\noindent $(ii)$ For each solution  $\trajxo(\cdot)$ of   \eqref{def:traj},  there exists a control
$\a(\cdot)=\big(\alpha_1(\cdot),\alpha_2(\cdot),\mu(\cdot)\big) \in \Ac$ such that
\begin{equation}\begin{aligned}\label{fond:traj}
\dottrajxo (t) &= b_1\big(\trajxo (t), \frac{\trajxo}{\e} (t),\alpha_1(t)\big)\mathds{1}_{\apg t  \:   : \:  \trajxo (t) \in  \e\Omega_1 \chg }+
b_2\big(\trajxo (t), \frac{\trajxo}{\e} (t),\alpha_2(t)\big)\mathds{1}_{\apg t  \:   : \:   \trajxo (t) \in \e \Omega_2 \chg } \\[2mm]
&+ b_\H\big(\trajxo (t) , \frac{\trajxo}{\e} (t),\a(t)\big)\mathds{1}_{\apg  t  \:   : \:  \trajxo (t) \in  \e\H \chg }\,, 
\end{aligned}\end{equation}
(where $\mathds{1}_{I}(\cdot)$ stands for the indicator function of the set $I$, and 
for the sake of simplicity the $\e$-dependence of the control $\a=\a^\e$ is not written.)

\noindent $(iii)$ Moreover, 
$$b_\H\big(\trajxo (t) ,\frac{\trajxo}{\e} (t),\a(t)\big)\cdot \nor_1( \frac{\trajxo}{\e} (t))= 0 \quad \hbox{a.e. on  }\{\trajxo(t) \in\e  \H \}\; .$$

\end{theo}

It is worth remarking that, in Theorem~\ref{def:dynAssume}, a solution $\trajxo(\cdot) $ can be associated to several controls $a(\cdot)$. Fix $\e >0$, to set properly the control problems, 
 we introduce the set $\Ta$ of admissible controlled trajectories starting from the initial datum $x_0$
\begin{equation*}
\Ta:= \big\{  (\trajxo(\cdot),\a(\cdot))\in {\rm Lip}(\R^+;\R^N) \times \A \mbox{ such that  }
\mbox{ \eqref{fond:traj}  is fulfilled and }  \trajxo(0)=x_0   \big\}
\end{equation*}
and we set (the $\e$ and $x_0$ dependence is not explicitly written)
$$
  \Euno:= \apg  t  \:   : \: \trajxo (t) \in \e\Omega_1   \chg,\quad \Edue:= \apg  t  \:   : \: \trajxo (t) \in \e\Omega_2  \chg,\quad
\Eh:= \apg  t  \:   : \: \trajxo (t) \in \e \H  \chg\,.
$$
We finally define the set of regular controlled trajectories
$$
\Treg:= \big\{  (\trajxo(\cdot),\a (\cdot)) \in \Ta \mbox{ such that, for almost all }  t\in\Eh,  \: b_\H(\trajxo(t),\frac{\trajxo}{\e} (t),a(t))   \mbox{ is regular}  \big\}.
$$
Recall that a \textit{regular} dynamics $b_\H(x,\y,\a)$ on $ \e
\H$ with $a=(\alpha_1, \alpha_2, \mu)$
is such that $b_i (x,\y,\alpha_1)\cdot  \nor_i(\y) \geq 0$ (where $\nor_i(\y)$ is the unitary normal exterior vector to $\Omega_i$ in $\y$) 
while a \textit{singular} a dynamic $b_\H(x,\y,\a)$ is such that
$b_i (x,\y,\alpha_1)\cdot  \nor_i(\y) < 0$.

%

\noindent {\bf The cost functional.}
Our aim is  to minimize an infinite horizon  cost functional such that we respectively pay $l_{i}$ if the trajectory is in $\e \Omega_i$, $i=1,2$ and
$l_\H$ if it is on $\e\H$. \\

More precisely,  the cost associated to $(\trajxo(\cdot) ,\a)  \in \Ta$ is
$$  
J(x_0; (\trajxo, \a)):=\int_0^{+\infty} l\big(\trajxo (t),\frac{\trajxo}{\e} (t),a(t) \big) e^{-\lambda t} dt
$$
where the Lagrangian is given by ($ l_\H$ is defined in \eqref{lH})
$$\begin{aligned}
    l(\trajxo (t),\frac{\trajxo}{\e} (t),a(t)) & :=  l_1(\trajxo (t),\frac{\trajxo}{\e} (t),\alpha_1(t)) \mathds{1}_{\Euno}(t)
       +  l_2(\trajxo (t),\frac{\trajxo}{\e} (t),\alpha_2(t)) \mathds{1}_{\Edue }(t) \\
       & +   l_\H(\trajxo (t),\frac{\trajxo}{\e} (t),\a(t))\mathds{1}_{\Eh }(t)\,.
\end{aligned}$$

For the sake of simplicity we also  set  ($ b_\H$ is defined in \eqref{bH})
$$\begin{aligned}
b(\trajxo (t),\frac{\trajxo}{\e}(t),a(t)) & :=  b_1(\trajxo (t),\frac{\trajxo}{\e} (t),\alpha_1(t)) \mathds{1}_{\Euno}(t)
       +  b_2(\trajxo (t),\frac{\trajxo}{\e} (t),\alpha_2(t)) \mathds{1}_{\Edue }(t) \\
       & +   b_\H(\trajxo (t),\frac{\trajxo}{\e} (t),\a(t))\mathds{1}_{\Eh }(t)\,.
\end{aligned}$$

\noindent {\bf The value functions. }For each initial data $x_0$, we define the following two value functions
\be \label{eq:valorem}
 \uem (x_0):= \inf_{(\trajxo,\a) \in \Ta} J(x_0; (\trajxo, \a))
 \ee
\be \label{eq:valorep}
 \uep (x_0):= \inf_{(\trajxo,\a) \in \Treg} J(x_0; (\trajxo, \a)).
 \ee


The most important consequence of the controllability assumption [H2]  is that both value functions  $\uem$ and   $\uep$ are uniformly Lipschitz continuous.

\begin{theo} \label{teo:unifLip}
Assume [H0], [H1], 
and [H2]. Then, the value functions $\uem$ and $\uep$  are  bounded, Lipschitz continuous functions
from $ \R^N$ into $\R$. Their $W^{1,\infty}$ norm is also uniformly
bounded with respect to $\e$.  
\end{theo}
{\bf Proof.}  For the details of the proof see Theorem  2.3 in \cite{BBC}. Here we only recall  that: if $M_l$ and $\delta$ are given in [H2] and [H3] 
then  the Lipschitz constant  is  $\frac{M_l}{\delta}$, hence it does not  depend on $\e$.
  \hfill $\Box$

The first key result is the {\bf Dynamic Programming Principle.}
\begin{theo}
Assume [H0], [H1], 
and [H2].
Let  $\uem,\uep$ be the  value functions defined in  \eqref{eq:valorem} and \eqref{eq:valorep}, respectively.
For each  initial data $x_0$, and each time $ \tau \geq 0$, we have
\begin{equation} \label{DPP-uem}
   \uem (x_0)= \inf_{  (\trajxo,\a) \in \Ta} \left\{ \int_0^\tau l\big(\trajxo (t),\frac{\trajxo}{\e} (t),a(t) \big)  e^{-\lambda t} dt
    +  e^{-\lambda \tau}  \uem (\trajxo(\tau)) \right\}
\end{equation}
\begin{equation} \label{DPP-uep}
   \uep (x_0)= \inf_{  (\trajxo,\a) \in \Treg} \left\{ \int_0^\tau  l\big(\trajxo (t),\frac{\trajxo}{\e} (t),a(t) \big)  e^{-\lambda t} dt
    +  e^{-\lambda \tau}  \uep (\trajxo(\tau)) \right\}.
\end{equation}
\end{theo}
{\bf Proof.} The proof is classical so we will omit it. \hfill $\Box$

As a consequence  of the DPP we  obtain that   both the value functions  $\uem$ and  $\uep$ are viscosity solutions of the {\bf Hamilton-Jacobi-Bellman} equation (\ref{Bellman-OmH}),
while  they  fulfill different  inequalities on the hyperplane $\H$.

\begin{theo}\label{teo:HJ}
Assume [H0], [H1], 
and [H2].
The value functions $\uem$ and   $\uep$ are  both  viscosity solutions of the Hamilton-Jacobi-Bellman equation \eqref{Bellman-OmH}. \\
Moreover, $\uem$ is a solution of  \eqref{Bellman-OmH-} and $\uep$ is a solution of  \eqref{Bellman-OmH+}.
%
\end{theo}
{\bf Proof.} The proof is given in \cite[Theorem 2.5]{BBC}  (see also \cite[Theorem 3.3]{BBC2}). \hfill $\Box$

\

We end this section by stating two comparison results we will need for \eqref{Bellman-OmH-}. 
The first one is a {\it strong comparison result} in $\R^N$,  while the second one is a {\it local comparison result} 
 we will need in the proof of the convergence result  Theorem \ref{Umeno-convergenza},  below. Moreover, we prove that 
 that  $\uem$ and $\uep$ are the minimal supersolution and the maximal subsolution of  \eqref{Bellman-OmH}, respectively.

\begin{theo}  \label{sous-max-sur-min} \label{Uni-dim-N} \label{Umeno-comparisonPalle}
Assume [H0], [H1], 
and [H2]. \\
(i) $\uem$ is the minimal supersolution and solution of \eqref{Bellman-OmH}.\\
(ii) $\uep$ is the maximal subsolution and solution of \eqref{Bellman-OmH}.\\
(iii) Let $u$ be a bounded, Lipschitz continuous subsolution of  \eqref{Bellman-OmH-} and $v$ be a bounded, lsc supersolution of  \eqref{Bellman-OmH-}. Then $u \leq v $ in $\R^N$.\\
(iv)  Fix $R>0$ and $\xi \in \R^N$. Let $u$ be a  bounded, Lipschitz continuous subsolution  and  $v$ a lsc supersolution of  \eqref{Bellman-OmH-}  for all $x \in B(\xi,R)$.
\be \label{tesi-comparisonPalle-um}
\mbox{ If }   u \leq  v \mbox{ on } \partial B(\xi,R)  \quad   \mbox{ then }  \quad    u \leq  v \mbox{ in }  B(\xi,R).
\ee 
(v)  Fix $R>0$ and $\xi \in \R^N$. Let $u$ be a  bounded, Lipschitz continuous subsolution of  \eqref{Bellman-OmH} for all $x \in B(\xi,R)$. 
\be \label{tesi-comparisonPalle-up}
\mbox{ If }   u \leq  \uep \mbox{ on } \partial B(\xi,R)  \quad   \mbox{ then }  \quad    u \leq  \uep \mbox{ in }  B(\xi,R).
\ee 
\end{theo}
{\bf Proof. }  In order to prove (i),(ii) and (iii) we remark that  the proof of  \cite[Theorem 4.1]{BBC}  (or \cite[Theorem 4.4]{BBC2}) is local, therefore it can be adapted easily to this case, so we will omit it. The local comparison (iv)  follows  directly from \cite[Theorem 4.1]{BBC2} while  (v) can be easily proved adapting  the  proof of   \cite[Theorem 4.4(ii)]{BBC2}.  In particular,  note that, since all the argument are local   we can use Theorem 3.5 and Theorem 3.7   in \cite{BBC2} to adapt the proof  of  \cite[Lemma 6.1]{BBC2}. (The same idea is detailed in  \cite[Theorem 4.1]{BBC}).
 \hfill $\Box$

\begin{rem}  \label{remsottosolHpiu}
We recall here that when we are dealing with $\mathbf{H^+}$ the two equations   \eqref{Bellman-OmH+}  and \eqref{Bellman-OmH}   are equivalent. 
Indeed in \cite[Theorem 3.7]{BBC2} is it proved that any subsolution  $u$ of    \eqref{Bellman-OmH}  fulfills 
$
\lambda u(x) + \HTreg(x,\frac{x}{\eps},Du)   \leq 0
$
 in the sense of Definition \ref{defi:sousolH}. However it is worth emphasizing the fact that $\HTreg$ plays a central role in the proof of Theorem 3.5 {\it (ii)} and  
 {\it (v)}, this is why we keep it here. 
\end{rem}

\section{ The homogenization result  for   $U_\eps^{-}$ }\label{sec:homoperumeno}

\subsection{The cell problem and the definition of the effective Hamiltonian.}
Le us first study the   {\bf cell problem}. 

\begin{theo}  \label{cell-Umeno-esistenza}
Assume [H0], [H1], 
and [H2].
For any $x,p\in \R^N$, there exist a unique constant $C=\bar{H}^-(x,p)$ such that the following cell problem has a Lipschitz continuous, $ \Z^N
$-periodic viscosity solution $\umcor$
\be \label{sys:1cellU-}
  \HHm(x,y,Dv(y)+p) =C     \hbox{ in   }  \R^N \,.\\
\ee
\end{theo}

To prove Theorem \ref{cell-Umeno-esistenza} we introduce the classical {\bf  $\rho$-problem}:
let $\rho >0$, for any $p\in \R^N$ and $x\in\R^N$ ($x$ and $p$ are the "frozen" variables), we denote by  $\ualpham$ the solution of

\be \label{sys-alpham}
  \rho u(y) + \HHm(x,y,Du(y)+p) = 0     \hbox{ in   }   \R^N \,.
\ee

As in Paragraph \ref{sec:ueps} (see also the results in \cite{BBC2}) it is possible to prove that there is one and only one solution $\ualpham$ of \eqref {sys-alpham}. Moreover  $\ualpham$ is characterized by  being the value function of an optimal control problem. 
For the sake of clarity we describe all these results in Section \ref{Appendice1} below.

\noindent{\bf Proof of Theorem \ref{cell-Umeno-esistenza}.}
The proof is classical. Thanks to the characterization \eqref{eq:valoralpham}, $\ualpham$ is a Lipschitz continuous, $ \Z^N
$-periodic function. Moreover, as a consequence of the uniform controllability condition in [H2],  the Lipschitz constant of $\ualpham$ can be chosen independent of $\rho$ (cf. the proof of Theorem \ref{teo:unifLip}).

Let us define  $\valpham(y)=\ualpham(y)-\ualpham(0)$ and $\lrm=-\rho \ualpham(0)$. By easy and classical estimates, the $\lrm$ are bounded and since the $\valpham$ are equi-Lipschitz continuous and periodic, they are also equi-bounded. Therefore, we can apply Ascoli-Arzel\`a Theorem and, up to extractions of subsequences, we may assume that $\{\valpham\}_\rho$ converges uniformly to a Lipschitz $\Z^N$-periodic fonction $\umcor$ and $\{\lrm\}_\rho$ converges to a constant $\bar{H}^-$.

Next we can use the stability result, Theorem \ref{teostab} and deduce that $\umcor$ and $\bar{H}^-$ satisfy \eqref{sys:1cellU-}.

To prove the uniqueness of the constant $\bar{H}^-$, we argue by contradiction, assuming that $(v,\nu)$  and $(w,\mu)$ are two different solutions of   \eqref{sys:1cellU-}.  
We suppose, for instance that $\nu < \mu$ and, without loss of generality, we can assume that $v >  w$ in $\R^N$ by adding a suitable large constant to $v$.
 
Let us fix $\rho$ small enough to have $ \rho v+\nu \leq  \rho w+\mu$  for all $x \in \R^N$.  
Since $(v,\nu)$ is a solution of \eqref{sys:1cellU-} we have that $v$ is a solution of 
\be 
 \rho v(y) + \HHm(x,y,Dv(y)+p) =  \rho v(y)  +\nu  \hbox{ in   }   \R^N
\ee
while $w$ is a solution of 
\be 
\rho w(y) + \HHm(x,y,Dw(y)+p) =   \rho w(y)+\mu      \hbox{ in   }   \R^N.
\ee
Thus, by  the comparison result, we have that $v \leq w $ in all $\R^N$ which is  a contradiction and the claim is proved.   
\hfill $\Box$

\subsection{Study of the $\rho$--control problem}\label{Appendice1}

In this section, we describe the control problem associated to the Bellman Equation~\eqref{sys-alpham} which is used in the proof of  Theorem  \ref{cell-Umeno-esistenza}.

As in Paragraph \ref{sec:ueps} 
in order to define the dynamic, we have to consider the solutions of  the differential inclusion 
\begin{equation}\label{def:trajyo}
\dottrajyoyo (t) \in \mB(x,\trajyoyo(t))  \quad \hbox{for a.e.  } t \in (0,+\infty)  \: ; \quad \trajyoyo (0)=y_0.
\end{equation}

We recall that
 there exist controls $\alpha_1(\cdot), \alpha_2(\cdot), a(\cdot)$ such that, for almost every $t \in \R$, one has
\begin{equation}\label{fond:trajyo}
\dottrajyoyo (t) =  b_1\big(x, \trajyoyo (t),\alpha_1(t)\big)\mathds{1}_{\Euno}(t)+
b_2\big(x, \trajyoyo (t),\alpha_2(t)\big)\mathds{1}_{\Edue}(t) 
+ b_\H\big(x ,\trajyoyo (t),\a(t)\big)\mathds{1}_{\Eh }(t)\,
\end{equation}
with (the dependence on $x$ and $y_0$ is not explicitely written)
$$
  \Euno:= \apg  t  \:   : \: \trajyoyo (t) \in \Omega_1   \chg,\quad \Edue:= \apg  t  \:   : \: \trajyoyo (t) \in \Omega_2  \chg,\quad
\Eh:= \apg  t  \:   : \: \trajyoyo (t) \in  \H  \chg\, .
$$
Finally we set 
\begin{equation*}
\Tayo : = \big\{  (\trajyoyo(\cdot),\a(\cdot))\in {\rm Lip}(\R^+;\R^N) \times \A \mbox{ such that  }
\mbox{ \eqref{fond:trajyo}  is fulfilled and }  \trajyoyo(0)=y_0   \big\}.
\end{equation*}

Following \cite{BBC2},  one can prove that there is one and only one solution $\ualpham$ of \eqref {sys-alpham} which is given by
\be \label{eq:valoralpham}
\ualpham (y_0):= \inf_{(\trajyoyo,\a) \in \Tayo} J^{\rho}(y_0; (\trajyoyo, \a))
 \ee
where
\be  \label{def-costo-y}
J^{\rho}(y_0; (\trajyoyo, \a)):=\int_0^{+\infty} \tilde{l}\big(x,p,\trajyoyo(t),a(t) \big) e^{-\rho t} dt
\ee
and
\begin{equation}\label{elletilde3}
\tilde{l}\big(x,p,\trajyoyo(t),a(t) \big)= l(x,\trajyoyo (t),a(t))+ b\big(x, \trajyoyo (t),a(t))\cdot p .
\end{equation}

We have the following result (see also \cite{BBC2}):
\begin{theo}   \label{teoconfnoto}
Assume [H0], [H1], 
and [H2]. Fix $x,p \in R^N$. For  $i=1,2$  let $f_i: \R^N \ds : \R$ be a bounded continuous function and $u_i$, $i=1,2$ be  the solution of the following equation :
\be 
 \rho u_i(y) + \HHm(x,y,Du_i+p) =   f_i(y)     \hbox{ in   }   \R^N,
\ee

$$ 
\mbox{if }f_1(x) \leq f_2(x) \mbox{ for any } x\in\R^N   \quad   \mbox{  then }  \quad  u_1(x)  \leq u_2(x)   \mbox{ for any  }x\in\R^N. 
$$
\end{theo}
{\bf Proof.}  The result follows from the following characterization of the solutions $u_i$:
 \begin{equation*} 
 u_i (y_0)= \inf_{(\trajyoyo,\a) \in \Tayo} \int_0^{+\infty} ( \tilde{l}\big(x,p,\trajyoyo(t),a \big) + f_i(\trajyoyo(t) )   ) e^{-\rho t} dt .
  \end{equation*}
\hfill $\Box$

\subsection{Properties of  $\bar{H}^-$ }

We complement this result by proving that  the effective Hamiltonian $\bar{H}^-$ fulfills the properties needed to  obtain a comparison result for the limit problem. 

\begin{theo} \label{propHeffe}
Assume [H0], [H1], 
and [H2]. 
Let $\bar{H}^- :\R^N \times \R^N \ds \R$ be defined in Theorem   \ref{cell-Umeno-esistenza}.  \\
(i)   There exists a modulus $w(\cdot)$ such that for any $x\in\R^N$, $z\in\R^N$,
\begin{equation}  \label{H-barre-lip:x}
\vert \bar{H}^- (x,p)-\bar{H}^- (z,p)\vert\leq w(x-z) (1+\vert p\vert).
\end{equation}
(ii) There exists a constant M such that for any $x\in\R^N$,  for any $p\in\R^N$, $q\in\R^N$
 \begin{equation}  \label{H-barre-lip:pH}
\vert  \bar{H}^- (x,p)- \bar{H}^- (x,q)\vert\leq M\vert p-q\vert.
\end{equation}
(iii) For any  $x\in\R^N$ and for any $p\in\R^N$ we have 
\[  \bar{H}^-(x,p)\geq -M_l+\delta\vert p\vert .  \]
\end{theo}
{\bf Proof. }The proofs are classical.
Proofs of (i) and (ii)  are based on the comparison principle for the $\rho$-problem defining $\bar{H}^- (x,p)$.\\
For instance, let us prove (i):\\
For  $\rho\in\R^+$, for any $p\in \R^N$, $x\in\R^N$ and $z\in\R^N$, we denote by  $\ualpham(x,\cdot,p)$ the solution of (\ref{sys-alpham})
and by $\ualpham(z,\cdot,p)$ the solution of
\be \label{sys-alpham-z}
 \rho u(y) + \HHm(z,y,Du(y)+p) = 0        \hbox{ in   }   \R^N.
\ee
Our aim is to prove that  $\ualpham(z,\cdot,p) +  \frac{1}{\rho}  w(x-z) (1+\vert p\vert)$ is a supersolution of  (\ref{sys-alpham}).  
Indeed, if this is true, using the comparison principle (Theorem  \ref{Uni-dim-N} (iii)) with $\ualpham(x,\cdot,p)$ (subsolution of  \eqref{sys-alpham}
satisfying  $\rho u(y) + \HT(x,y,D_\H u+p)  \leq 0$) we deduce that
\[
\ualpham(x,y,p)\leq \ualpham(z,y,p)+\frac{1}{\rho} w(x-z) (1+\vert p\vert).
\]
Multiplying by $\rho$ and letting $\rho \rightarrow 0$, we have 
\[
\bar{H}^- (x,p)\geq \bar{H}^- (z,p)- w(x-z) (1+\vert p\vert),
\]
and, reversing the roles of $x$ and $z$,  we conclude the proof.  \\
We  prove  now that $\ualpham(\cdot, z,p) +  \frac{1}{\rho}  w(x-z) (1+\vert p\vert)$ is a supersolution of  (\ref{sys-alpham}).  \\
Since the argument is completely similar we detail only the case $y\in\Omega_i$. 

Since  $\ualpham(z,\cdot,p)$ is a solution of \eqref{sys-alpham-z}, by coerciveness  of $H_i$ there exists three constants $K_1>0$, $K_2$, $K_3$ 
(depending only on the constant defined in [H0]...[H2]) such that 
$$
K_2+K_1 \vert D_y\ualpham(z,y,p) +p\vert \leq  \vert  - \rho \ualpham(z,y,p) \vert. 
$$
and
$$ 
\vert \ualpham(z,y,p)  \vert 
\leq\frac {K_3}{\rho} (1 + \vert p \vert). 
$$
We deduce that  there exists a constant $K$ such that 
\begin{equation}
\label{eq:K123}
\vert D_y\ualpham(z,y,p) +p\vert \leq K (1+ \vert p \vert).
 \end{equation}

\noindent Moreover there exist a modulus $w$ such that for any $q\in \R^N$, $x\in  \R^N$, $z\in  \R^N$,
\[ H_i(z,y,q) \leq H_i(x,y,q) +(1+ \vert q\vert)w(x-z).\]
Therefore, in the viscosity sense:
\begin{equation}
\begin{split}   &\rho \ualpham(z,y,p)+ H_i(x,y,D_y\ualpham(z,y,p)+p) \\ 
&\geq\rho \ualpham(z,y,p)+ H_i(z,y,D_y\ualpham(z,y,p)+p)-(1+\vert D_y\ualpham(z,y,p) +p\vert)w(x-z)\\
&  \geq-(1+\vert D_y\ualpham(z,y,p) +p\vert)w(x-z)
 \end{split}
 \end{equation}
Thanks to \eqref{eq:K123} we  have  that there exist a modulus $w$ (note that  this is not exactly the same as before) such that
  $ \ualpham(z,y,p)+\frac{1}{\rho}(1+\vert p\vert)w(x-z)$ satisfies 
\[\rho v(y) + H_i(x,y,D_yv(y)+p)\geq0, \hbox{ in }\Omega_i
\]in the viscosity sense and (i) holds.   \\
The coerciveness property (iii) follows from
\[
-\rho \ualpham(x,y,p)\geq -M_l+\delta\vert p\vert
\] 
letting $\rho\rightarrow 0$.
Indeed,
thanks to periodicity and continuity, $\ualpham(x,y,p)$ has a maximum point and a minimum point called respectively $y_M$ and $y_m$.\\
Then  $-\rho \ualpham(x,y_m,p)\geq -\rho \ualpham(x,y,p)\geq -\rho \ualpham(x,y_M,p)$. \\
Using only Definition  \eqref{progenNOH} of sub and super-solution, there exist $i\in \{ 1,2\}$ and $j\in \{ 1,2\}$ such that taking $\phi=0$ as test function in the equations satisfied by $V_\rho^-$ we have
\[
-\rho \ualpham(x,y_M,p)\geq H_i(x,y_M,p), \quad \mbox{ and }  
 \quad H_j(x,y_m,p)\geq -\rho \ualpham(x,y_m,p)
\]
this implies
\[
H_j(x,y_m,p)\geq -\rho V_\rho^-(y) \geq H_i(x,y_M,p)
\]
and thanks to the controllability condition [H2] and the boundedness of $l_i$ in [H1] it is easy to remark that
\[
H_i(x,y_M,p)\geq -M_l+\delta\vert p\vert.
\]

\hfill $\Box$

\subsection{The convergence result}

Before proving the convergence result we state here a comparison result for the limiting problem.

\begin{theo}  \label{Umeno-unicit? limite}
Assume [H0], [H1], 
and [H2].  Let $\bar{H}^-$ be defined in Theorem   \ref{cell-Umeno-esistenza}. Let  $u$ and $v$ be respectively a  bounded usc subsolution and a bounded lsc supersolution of   
\be \label{sys:limiteUmeno }
  \begin{array}{cc}
{\dis  \lambda w(x) + \bar{H}^-(x,Dw(x))  =0   }&   \hbox{ in   }  \R^N \,, \\
\end{array}  
\ee
Then $u(x) \leq v(x)$ in  $\R^N. $
\end{theo}

\noindent{\bf Proof of Theorem \ref{Umeno-unicit? limite}.}   
Thanks to Proposition   \ref{propHeffe} we are able to apply the classical  comparison results for Hamilton-Jacobi-Bellmann Equation in $\R^N$.
(See for instance \cite{BCD},  \cite{Ba},  \cite{BaCIME}, \cite{HICIME} or \cite{L}).
\hfill $\Box$  

\

We are finally ready to prove the  convergence result. More precisely 
\begin{theo}  \label{Umeno-convergenza}
Assume [H0], [H1] 
and [H2]. Let $\bar{H}^- :\R^N \times \R^N \ds \R$ be defined in Theorem   \ref{cell-Umeno-esistenza}. The sequence  $(\uem)_{\e >0}$ converges locally uniformly in $\R^N$ to a function $U^-$ which is the unique solution of 
\be \label{sys:limiteUmeno}
  \begin{array}{cc}
{\dis  \lambda U^-(x) + \bar{H}^-(x,DU^-)  =0   }&   \hbox{ in   }  \R^N \,.\\
\end{array} 
\ee
 \end{theo}
\noindent{\bf Proof
.}   We first remark that, in view of  Theorem \ref{teo:unifLip}, the functions $\uem$ are bounded and Lipschitz continuous uniformly with respect to $\e$.  Therefore, thanks to Ascoli-Arzela's Theorem, we may assume, up to the extraction of a subsequence, that the  sequence $(\uem)_{\e >0}$ converges locally uniformly in $\R^N$ to a bounded, Lipschitz continuous function $U^-$. Because of Theorem~\ref{Umeno-unicit? limite}, in order to conclude, we only need  to prove that $U^-$ is a sub and a supersolution of 
 \be 
  \begin{array}{cc}
{\dis  \lambda U^-(x) + \bar{H}^-(x,DU^-)  =0   }&   \hbox{ in   }  \Omega \,.\\
\end{array} 
\ee
Since both proofs use the same technique we will detail only the proof of  $U^-$  being  a supersolution.

Let  $\phi$ be  a $C^1$-function in $\R^N$  and  $\bar{x}$ be a local strict minimum point of $U^--\phi$ : we may assume without loss of generality that
\be  \label{ipminstrettoSuper-um}
\hbox{there exists  } \bar{r} >0  \hbox{ such that } (U^--\phi) (x) > 0 \; \hbox{for all  }x \in B(\bar{x},\bar{r}) \; \hbox{and} \; (U^--\phi) (\bar{x})=0.
 \ee
Our aim is to prove that 
 \be \label{tesi-UmSuper}
  \lambda  \phi(\bar{x})+  \bar{H}^-(\bar{x},D \phi (\bar{x}))    \geq 0.
 \ee
We argue by contradiction, assuming that  
  \be \label{tesi-Um-assurdoSuper}
  \lambda  \phi(\bar{x})+  \bar{H}^-(\bar{x},D \phi (\bar{x}))    \leq -\eta  < 0.
 \ee

Let $\Hbarm(\bar{x}, D\phi(\bar{x}))$  be defined  as in  Theorem \ref{cell-Umeno-esistenza} and $\umcor$ be the related unique solution of  system \eqref{sys:1cellU-}, i.e.
\be  \label{sysV-}
\HHm(\bar{x},y,D\umcor (y)+  D \phi (\bar{x}))=\bar{H}^-(\bar{x},D \phi (\bar{x}))  \quad \mbox{ in } \R^N.
\ee
 
\begin{lem}\label{keypointconv} There exist $\e_0>0$,  $r_0>0, \gamma_0>0$  such that the function $\chi_\e(x):=\phi(x)+ \e \umcor( \frac{x}{\e}) + \gamma$ is a viscosity subsolution of   \eqref{Bellman-OmH-} for all $x \in B(\bar{x},r)$, for any $r\leq r_0$ and $\gamma\leq \gamma_0$,  if $\e \leq \e_0$.
\end{lem} 

\noindent{\bf Proof of Lemma~\ref{keypointconv} :} 
Fix $\e,r >0$ and a point $x \in B(\bar{x},r)$. Since  $\umcor(\cdot)$ is a viscosity solution of  \eqref{sysV-} we have that 
$$  
\HHm(\bar{x},\y,D\umcor \big(\y \big)+  D \phi (\bar{x}))=\bar{H}^-(\bar{x},D \phi (\bar{x}))  
$$
in the viscosity sense. We remark now that, 
$$ 
\HHm(\bar{x},\y,D\umcor \big(\y \big)+  D \phi (\bar{x}))=\HHm(x,\y,D\chi_\e(x)) + O(r)
$$
thanks to the Lipschitz properties of $b_i$ and $l_i$ for  $H_i$ and Proposition  \ref{propHtan}  for $\H_T$. 
Moreover, since  $\umcor$ is bounded and  $\phi$ is regular we have 
$
 \lambda \chi_\e(x)= \lambda \gamma + O(\e) +O (r) + \phi(\bar{x}). 
$
Therefore, by assumption  \eqref{tesi-Um-assurdoSuper} we can deduce the following inequality  
 in the viscosity sense  
\begin{eqnarray*}
 \lambda \chi_\e(x) + \HHm(x, \frac{x}{\e},D \chi_\e(x))&=&    
  \lambda \gamma + O(\e) +O (r) + \phi(\bar{x}) + \HHm(\bar{x},\y,D\umcor \big(\y \big)+  D \phi (\bar{x})) \\
  &=&   \lambda \gamma + O(\e) +O (r) + \phi(\bar{x}) +\bar{H}^-(\bar{x},D \phi (\bar{x}))   \\
  & \leq & \lambda \gamma  +    O(\e)+O(r) - \eta \leq   O(\e)-\frac{\eta}2\; ,
\end{eqnarray*}
if, say, $\lambda \gamma \leq \eta/4$, for $r$ and $\e$ small enough, depending only on $\eta$.  The proof of the Lemma is then completed.   \hfill $\Box$

Now we note that, by  \eqref{ipminstrettoSuper-um}, if $r\leq \bar{r}$, then there exists $\gamma_r >0$, such that  
$$
  \lim_{\e \ds 0} \uem(x) = U^-(x) \geq \phi(x) +\gamma_r   \quad \hbox{for all   }x \in  \partial B(\bar{x}, r)\; .
$$
Since the limit is uniform and $\umcor( \cdot)$ is bounded on $\R^N$,  for $\eps$ small enough, we have
$$
\uem(x) \geq \phi(x) +\gamma_r  \geq \phi(x) + \eps \umcor( \frac{x}{\eps})  +\frac{\gamma_r}{2} \quad\hbox{for all   }x \in  \partial B(\bar{x}, r) \; .
  $$ 
 
Choosing now $\gamma, r$ small enough ($\gamma\leq\gamma_r /2$ and $r$ small enough), we can have, at the same time, 
\begin{itemize}
\item[(i)] $\chi_\e(x) $ is a subsolution of  
\eqref{Bellman-OmH-} in $B(\bar{x},r)$ .
 \item[(ii)] $\chi_\e(x) $ is less than $\uem$ on $\partial B(\bar{x},r)$ .
\end{itemize}
%

Therefore, applying the comparison result  (cf. Theorem~\ref{Umeno-comparisonPalle}), we obtain
$$
\uem(x)   \geq \phi(x) + \eps \umcor( \frac{x}{\eps})  + \frac{\gamma}{2} \quad  \hbox{for all  } x \in   B(\bar{x}, \bar{r})  \; ,
  $$  
and letting $\eps \ds 0$  we have 
$$
U^-(x) \geq  \phi(x) + \frac{\gamma}{2} \quad  \hbox{for all  } x \in   B(\bar{x}, \bar{r})  \; .
$$
Taking $x =\bar{x}$ we get $U^-(\bar{x}) \geq \phi(\bar{x}) + \frac{\gamma}{2}$ which is in contradiction with \eqref{ipminstrettoSuper-um} and the proof is complete.

\hfill $\Box$


\section{ The homogenization result  for   $U_\eps^{+}$  }\label{sec:homoperupiu}

As remarked in the Introduction, the  function $U_\eps^{+}$  can only be characterized by being  the maximal solution of \eqref{Bellman-Om}-\eqref{Bellman-H-sub}-\eqref{Bellman-H-sup}. Therefore we are going to prove the homogenization result  by performing all the arguments on the control formulas. However, since we are closely following the ideas of pde argument's we start again by studying the cell problem.

\subsection{The cell problem}

The following theorem plays the same role as Theorem \ref{cell-Umeno-esistenza} for $\Hbarm(x,p)\in \R$:
\begin{theo}  \label{cell-Upiu-esistenza}
Assume [H0], [H1], 
and [H2]. For any $x, p\in \R^N$, there exists a unique constant $\Hbarp(x,p)\in \R$ such that there exists a Lipschitz continuous, periodic function $V^+$ satisfying, for any $\tau\geq 0$ and $y_0\in\R^N$
\begin{equation} \label{H+v} 
V^+(y_0)=
 \inf_{(\trajyoyo,\a) \in \Tayoreg}  \left\{ \int_0^\tau \left(\tilde{l}\big(x,p,\trajyoyo(t),a(t) \big)
    +\Hbarp(x,p) \right)dt+V^+(\trajyoyo(\tau))
    \right\}
\end{equation}
where $\tilde{l}$ is defined in (\ref{elletilde3}). Moreover $V^+$ is a viscosity subsolution of
\be \label{sys:1cellU+ }
   \HHp(x,y,DV^++p) =\Hbarp(x,p)      \hbox{ in   }  \R^N.  \footnote{Note that, thanks to Remark  \ref{remsottosolHpiu} we only have to prove that $V^+$  is a subsolution of   $H(x,y,DV^++p) =\Hbarp(x,p)$.}
\ee
 Finally, for all $y_0\in\R^N$ we have
\begin{equation}  \label{formulasinto}
 \Hbarp(x,p)= \lim_{t \rightarrow +\infty}  \apt -  \inf_{(\trajyoyo,\a) \in \Tayoreg}  \big\{  \frac{1}{t}     \int_0^t \tilde{l}\big(x,p,\trajyoyo(t),a(t) \big)  dt  \big\}
     \cht.
\end{equation} 
 \end{theo} 
\begin{rem}   \label{remformrap}
Note that, since 
$$\tilde{l}\big(x,p,\trajyoyo(t),a(t) \big) =l\big(x,\trajyoyo(t),a(t) \big)  + b\big(x,\trajyoyo(t),a(t) \big)  \cdot p=l\big(x,\trajyoyo(t),a(t) \big) 
+\dot{Y}_{y_0}(t) \cdot p,  $$  
formula  \eqref{formulasinto}  can be rewritten as 
 \begin{equation}  \label{formulasintosempli}
 \Hbarp(x,p)= \lim_{t \rightarrow +\infty}  \apt -  \inf_{(\trajyoyo,\a) \in \Tayoreg}  \big\{  \frac{1}{t}     \int_0^t  l\big(x,\trajyoyo(t),a(t) \big)  dt
   +\frac{(\trajyoyo(t)-y_0)}{t}  \cdot p   \big\} \cht .
\end{equation}
Moreover, by taking the infimum on the  set of all  the trajectories  $\Tayo$(instead of 
 the set of  regular trajectories  $\Tayoreg$) we can obtain the  same characterization  for  $\Hbarm(x,p)$.  
\end{rem}

{\bf Proof}
We introduce the classical  $\rho$-problem.
Let $\rho\in\R^+$, for any fixed $p\in \R^N$ and $x\in\R^N$ 
\be \label{eq:valoralphap}
\Vrho (y_0):= \inf_{(\trajyoyo,\a) \in \Tayoreg} J^{\rho}(y_0,(\trajyoyo, \a))
 \ee
where $y_0\in\R^N$ and the cost function is defined as in \eqref{def-costo-y} but here the infimum is taken considering only the {\it regular} trajectories 
\begin{equation*}
\Tayoreg := \big\{  (\trajyoyo(\cdot),\a (\cdot)) \in \Tayo \mbox{ such that, for almost all }  t\in\Eh,  \: b_\H(x,\trajyoyo(t),a(t))   \mbox{ is regular} \big\}.
\end{equation*}
(Recall that a control $a$ is regular if   
$b_i (x,y,\alpha_i)\cdot  \nor_i(y) \geq 0$, $i=1,2$.)

By Theorem \ref{teo:HJ} and Theorem  \ref{sous-max-sur-min} we know that  $\Vrho$  is the  maximal 
subsolution of
\be \label{sys-alphaP}
  \rho u(y) + H^+(x,y,Du(y)+p) = 0       \hbox{ in   }   \R^N .
\ee

Thanks to Definition \eqref{eq:valoralphap}, $\Vrho$ is a $ \Z^N
$-periodic function and since the $b_i,l_i$ ($i=1,2$) are bounded, then an easy estimate proves that $\rho \Vrho$ is bounded, uniformly in $\rho$. Moreover, thanks to the uniform controllability condition in [H2], as proved in Theorem \ref{teo:unifLip}, $\Vrho$ is Lipschitz continuous and its Lipschitz constant is independent of $\rho$.

Let us define $\valphap(y)=\Vrho(y)-\Vrho(0)$ and $\walphap=-\rho \Vrho(0)$. Up to an extraction of a subsequence, thanks to Ascoli-Arzel\`a Theorem, we may assume that $(\valphap(y))_\rho$ converges uniformly in $\R^N$
to a Lipschitz $ \Z^N
$-periodic fonction $V^+$ and $(\walphap)_\rho$ converges to a constant that we will denote by  $\bar{H}^+(x,p)$.

Using the stability property (cf. Theorem \ref{teostabp}), we have  that $V^+$ is a subsolution of \eqref{sys:1cellU+ }. Moreover, by the Dynamic Programming Principle for $\Vrho$, we have, for each $y_0\in \R^N$ and each time $ \tau \geq 0$, 
\begin{equation} \label{DPP+}
\Vrho (y_0)= \inf_{(\trajyoyo,\a) \in \Tayoreg}  \left\{ \int_0^\tau  \tilde{l}\big(x,p,\trajyoyo(t),a(t) \big) e^{-\rho t} dt
    +  e^{-\rho\tau}  \Vrho (\trajyoyo(\tau)) \right\} \; .
\end{equation}
Next we use that
$$ \Vrho (0) = \int_0^\tau \rho\Vrho (0)e^{-\rho t} dt + \Vrho (0)e^{-\rho \tau }\,$$
which yields, by subtracting
$$
\valphap (y_0)= \inf_{(\trajyoyo,\a) \in \Tayoreg}  \left\{ \int_0^\tau  \left[\tilde{l}\big(x,p,\trajyoyo(t),a(t) \big) -\rho\Vrho (0)\right]e^{-\rho t} dt
    +  e^{-\rho\tau}  \valphap  (\trajyoyo(\tau)) \right\}\; .
$$
Using the uniform convergence of the sequence $\{\valphap\}_\rho$, it is easy to pass to the limit in this equality and to get
\begin{equation}  
V^+(y_0)=
 \inf_{(\trajyoyo,\a) \in \Tayoreg}  \left\{ \int_0^\tau \left( \tilde{l}\big(x,p,\trajyoyo(t),a(t) \big)
    +\bar{H}^+(x,p)\right)dt+V^+ (\trajyoyo(\tau))
    \right\}
\end{equation}
It is worth pointing out that here $\tau$ is arbitrary thus  $V^+$ satisfies this property for any $\tau \geq 0$.

Now we want to prove the uniqueness of $\bar{H}^+(x,p)$.


Let us suppose now that there exist two constants $\bar{H}^{+,1}(x,p)$ and $\bar{H}^{+,2}(x,p)$ such that there exist respectively two continuous periodic functions $v_1$ and $v_2$
safisfying \eqref{H+v}, then using that $\inf(\cdots)-\inf(\cdots) \leq \sup(\cdots)$, we obtain 
$$
v_1(y_0)-v_2(y_0)\leq  \sup_{(\trajyoyo,\a) \in \Tayoreg}  \left\{ \int_0^\tau(\bar{H}^{+,1}(x,p)-\bar{H}^{+,2}(x,p))dt + (v_1-v_2)(\trajyoyo(\tau)) \right\}\; .
$$
Thanks to the properties of periodicity and continuity of $v_1$ and $v_2$ there exist a $y^*$ such that 
$$
\max_{y\in \R^N}(v_1-v_2)(y)=v_1(y^*)-v_2(y^*)\; ,
$$
and we can use the preceding equality with $y_0=y^*$.

This leads to 
$$
v_1(y^*)-v_2(y^*)\leq  (\bar{H}^{+,1}(x,p)-\bar{H}^{+,2}(x,p))\tau + \max_{y\in \R^N}(v_1-v_2)(y)\; ,$$
and for this inequality to hold for $\tau >0$, this clearly implies that $\bar{H}^{+,1}(x,p)-\bar{H}^{+,2}(x,p) \geq 0$.
Exchanging the roles of $v_1$ and $v_2$, we obtain the opposite inequality, i.e. $\bar{H}^{+,1}(x,p)=\bar{H}^{+,2}(x,p)$.

We end the proof by remarking that we can deduce   \eqref{formulasinto} by  \eqref{H+v}  thanks to the boundness  of   $V^+$.  
\hfill $\Box$

\begin{rem}We remark that, not only  the above proof just requires that there exists a positive $\tau$ such that  \eqref{H+v} holds to obtain the uniqueness of $\bar{H}^+(x,p)$ but one can also use this proof to obtain further results. For example, if $\tilde w$ is a subsolution of (\ref{sys:1cellU+ }) associated to $\tilde{H}^{+}$, one can prove that it satisfies a suboptimality principle, i.e.
\begin{equation} \label{H+w} \tilde w(y_0)\leq
 \inf_{(\trajyoyo,\a) \in \Tayoreg}  \left\{ \int_0^\tau \left( \tilde{l}\big(x,p,\trajyoyo(t),a(t) \big)
    +\tilde{H}^{+}(x,p)\right)dt+\tilde w(\trajyoyo(\tau)).
    \right\}
    \end{equation}
Moreover, since $\tilde w$ is bounded the above argument with $v_1=\tilde w$ and $v_2=V^+$ leads to $\tilde{H}^{+}(x,p)\geq\bar{H}^{+}(x,p)$. This means that $\bar{H}^{+}=\inf \tilde{H}^{+}$,
    where the infimum is taken on the set of  subsolutions of (\ref{sys:1cellU+ }).
\end{rem}

Now we provide additional properties for $\bar{H}^{+}(x,t)$ to state that the effective Hamiltonian  fulfills the properties needed to  obtain a comparison result for the limiting problem
as in Theorem \ref{propHeffe} for $\bar{H}^-$.
\begin{theo}  \label{Car-Hbarra+0}Assume [H0], [H1], 
and [H2]. \\
(i) There exists a modulus $w(\cdot)$ such that for any $x\in\R^N$, $z\in\R^N$,
\begin{equation}  \label{H-barre+lip:x}
\vert \bar{H}^+ (x,p)-\bar{H}^+ (z,p)\vert\leq w(x-z) (1+\vert p\vert).
\end{equation}
(ii) For any $x\in\R^N$,  for any $p\in\R^N$, $q\in\R^N$
 \begin{equation}  \label{H-barre+lip:pH}
\vert  \bar{H}^+ (x,p)- \bar{H}^+ (x,q)\vert\leq M_b\vert p-q\vert.
\end{equation}
(iii) For any  $x\in\R^N$ and for any $p\in\R^N$ we have 
\[
\bar{H}^+(x,p)\geq -M_l+\delta\vert p\vert .
\]
\end{theo}

The proof is similar to that of Theorem \ref{propHeffe} but we can't use a comparison principle and we replace it by the use of maximal subsolutions.\\
{\bf Proof.}
For instance, let us prove (i):\\
Let us remember that $\Vrho$, that here we will denote by 
$\Vrho(x,\cdot,p)$, is the  maximal 
subsolution of \eqref{sys-alphaP}:
\be 
  \rho u(y) + H^+(x,y,Du(y)+p) = 0       \hbox{ in   }   \R^N
\ee
while 
$\Vrho(z,\cdot,p)$  is the  maximal 
subsolution of \label{sys-alphaP-z}
\be 
 \label{sys-alphaP-z}
  \rho u(y) + H^+(z,y,Du(y)+p) = 0       \hbox{ in   }   \R^N.
\ee
We want to prove that there exists a modulus $w$ such that the function 
$$v(y):=\Vrho(x,y,p)-\frac{w(x-z)}{\rho}(1+\vert p \vert)$$
is a subsolution of \eqref{sys-alphaP-z} (in the sense of Definition \ref{solviscogenNOH}). If this is true $v(y)\leq \Vrho(z,y,p)$, this implies that
$$\rho(\Vrho(x,y,p)-\Vrho(z,y,p))\leq  w(x-z)(1+\vert p \vert)$$
 letting $\rho\rightarrow 0$ we have $\Hbarp(z,p)-\Hbarp(x,p)\leq w(x-z)(1+\vert p \vert)$ and the proof is completed.\\
Let us suppose that $y\in\H$ (the other cases are analogous and even simpler), there exists $i=1$ (or $2$)  such that in the sense of viscosity:
\be \label{pippo} 
  \rho \Vrho(x,y,p) +H_1(x,y,D\Vrho(x,y,p)+p) \leq  0. 
\ee
Thanks to the coerciveness properties of $H_1$ there exist a constant $K_1>0$ and a constant $K_2$ such that $H_1(x,y,D\Vrho(x,y,p)+p)\geq K_1\vert D\Vrho(x,y,p)+p \vert +K_2$
and using the definition of $\Vrho(x,y,p)$ there exists a constant $K_3$ such that $\vert \rho \Vrho(x,y,p)\vert\leq K_3(1+\vert p \vert)$
\begin{equation*}
K_1\vert D\Vrho(x,y,p)+p \vert +K_2 \leq H_1(x,y,D\Vrho(x,y,p)+p) \leq  -  \rho \Vrho(x,y,p) \leq K_3(1+\vert p \vert)
\end{equation*}
then there exists a constant $K_4$ depending only on the constants defined in [H0], [H1] and [H2] such that
\begin{equation*}
\vert D\Vrho(x,y,p)+p \vert  \leq K_4(1+\vert p \vert).
\end{equation*}

This implies, using the regularity properties of $H_1$ and \eqref{pippo}  that there exists a modulus $w$ such that
\begin{equation*}
  \rho \Vrho(x,y,p) + H_1(z,y,D\Vrho(x,y,p)+p)\leq  H_1(z,y,D\Vrho(x,y,p)+p) - H_1(x,y,D\Vrho(x,y,p)+p)
\end{equation*}
 \begin{equation*}
\leq  C(1+\vert D\Vrho(x,y,p)+p\vert) w( x-z)\leq (1+\vert p \vert) w( x-z)
\end{equation*}
therefore the function $\Vrho(x,y,p)- \frac {  w( x-z)}{\rho} (1+\vert p \vert)$ verify in the viscosity sense
\begin{equation*}
  \rho u(y) + \min\{H_1(z,y,Du(y)+p),H_2(z,y,Du(y)+p)\}\leq  0.
\end{equation*}
\hfill $\Box$

\subsection{The convergence result}
We are now ready to prove the  convergence result. More precisely
\begin{theo}  \label{Upiu-convergenza}
Assume [H0], [H1], 
and [H2]. Let $\Hbarp$  be defined  as in  Theorem   \ref{cell-Upiu-esistenza} and $\uep$ as in \eqref{eq:valorep}. Then sequence  $(\uep)_{\e >0}$ converges locally uniformly in $\R^N$ to a continuous function $\up$, which is the unique viscosity solution of
\be
 \label{sys:limiteUpiu}
\dis  \lambda u(x) + \bar{H}^+(x,Du(x))  =0      \hbox{ in   }  \R^N .
\ee
 \end{theo}

{\bf Proof}  We first remark that, in view of  Theorem \ref{teo:unifLip}, the functions   $\uep$   are  equi-bounded and equi- Lipschitz continuous.  Therefore, by Ascoli-Arzela Theorem, we can extract a subsequence, still denoted by $(\uep)_{\e >0}$, which converges locally uniformly  on $\R^N$ to a function $\up$. 
 Our aim is then to prove that $\up$ is a solution of  \eqref{sys:limiteUpiu}. If this is the case, since by Theorem~\ref{Car-Hbarra+0}, \eqref{sys:limiteUpiu} has a unique viscosity solution (because $\bar{H}^+$ satisfies the classical assumption of classical comparison results)  then the whole sequence will converge to $\up$.
 
Since some parts of the proof are rather technical, we split it into three steps. The first step, concerning the supersolution property for  $\up$ is rather similar to the analogous proof for $\um$,
the principal tool is the local comparison principle stated in Theorem~\ref{Umeno-comparisonPalle} (v), unfortunately  to prove the subsolution property for  $\up$ will be more difficult because  a comparison principle concerning supersolutions and $\uep$ does not hold, in the second step we prove the subsolution property in the case when $b$ and $l$ do not depend on the first variable and in the last step  the subsolution property is proved in the general case 
using a sequence of approximating problems. 
 
 \bigskip
 
\noindent {\bf Step 1  : $\up$  is a supersolution of \eqref{sys:limiteUpiu}.} In this step, we follow readily the pde arguments which are already used for $\um$.
 
Let  $\phi$ be  a $C^1$-function in $\R^N$  and  $\bar{x}$ be a local strict minimum point of $\up-\phi$ : we may assume without loss of generality that
\be  \label{ipminstrettoSuper-up}
\hbox{there exists  } \bar{r} >0  \hbox{ such that } (\up-\phi) (x) > 0 \; \hbox{for all  }x \in B(\bar{x},\bar{r}) \; \hbox{and} \; (\up-\phi) (\bar{x})=0.
 \ee
Our aim is to prove that 
 \be \label{tesi-UpSuper}
  \lambda  \phi(\bar{x})+  \bar{H}^+(\bar{x},D \phi (\bar{x}))    \geq 0
 \ee
where  $\Hbarp(\bar{x}, D\phi(\bar{x}))$  is defined  as in  Theorem   \ref{cell-Upiu-esistenza}.
We argue by contradiction, assuming that  
  \be \label{tesi-Up-assurdoSuper}
  \lambda  \phi(\bar{x})+  \bar{H}^+(\bar{x},D \phi (\bar{x}))\leq -\eta  < 0.
 \ee
Let $\upcor(y)$ be the  subsolution  related  to  $\Hbarp(\bar{x}, D\phi(\bar{x}))$ fulfilling  \eqref{sys:1cellU+ },
  i.e.
\be  \label{sysV+}
\HHp(\bar{x},y,D\upcor (y)+  D \phi (\bar{x}))=\bar{H}^+(\bar{x},D \phi (\bar{x}))  \quad \mbox{ in } \R^N.
\ee
Let us first prove the following Lemma.
\begin{lem}\label{keypointconvp} There exists $r, \gamma>0$ small enough such that the function $\chi_\e(x):=\phi(x)+ \e \upcor( \frac{x}{\e}) + \gamma$ is a  subsolution of  \eqref{Bellman-OmH+} for all $x \in B(\bar{x},r)$, if $\e >0$ is small enough.
\end{lem} 

\noindent{\bf Proof of Lemma~\ref{keypointconvp} :} 
Fix $\e,r >0$ and a point $x \in B(\bar{x},r)$. Since  $\upcor(\cdot)$ is a  subsolution of  \eqref{sysV-} we have that 
$$  
\HHp(\bar{x},\y,D\upcor \big(\y \big)+  D \phi (\bar{x})) \leq \bar{H}^+(\bar{x},D \phi (\bar{x}))  
$$
in the viscosity sense. We remark now that, 
$$ 
\HHp(\bar{x},\y,D\upcor \big(\y \big)+  D \phi (\bar{x}))=\HHp(x,\y,D\chi_\e(x)) + O(r)
$$
thanks to the Lipschitz properties of $b_i$ and $l_i$ for  $H_i$ and  Remark \ref{propHtanreg}  for $\HTreg$. 
Moreover, since  $\upcor$ is bounded and  $\phi$ is regular we have 
$
 \lambda \chi_\e(x)= \lambda \gamma + O(\e) +O (r) + \phi(\bar{x}). 
$
Therefore, by assumption  \eqref{tesi-Up-assurdoSuper} we can deduce the following inequality  
 in the viscosity sense  
\begin{eqnarray*}
 \lambda \chi_\e(x) + \HHp(x, \frac{x}{\e},D \chi_\e(x))&=&    
  \lambda \gamma + O(\e) +O (r) + \phi(\bar{x}) + \HHp(\bar{x},\y,D\upcor \big(\y \big)+  D \phi (\bar{x})) \\
  &\leq &   \lambda \gamma + O(\e) +O (r) + \phi(\bar{x}) +\bar{H}^+(\bar{x},D \phi (\bar{x}))   \\
  & \leq & \lambda \gamma  +    O(\e)+O(r) - \eta \leq  0\; ,
\end{eqnarray*}
if, say, $\lambda \gamma \leq \eta/2$, for $r$ and $\e$ small enough, depending only on $\eta$.  The proof of the Lemma is then completed.   \hfill $\Box$

\

Now we remark that, by  \eqref{ipminstrettoSuper-up}, if $r\leq \bar{r}$, then there exists $\gamma_r >0$, such that  
$$
  \lim_{\e \ds 0} \uep(x) = U^+(x) \geq \phi(x) +\gamma_r   \quad \hbox{for all   }x \in  \partial B(\bar{x}, r)\; .
$$
Since this limit is uniform and $\upcor( \cdot)$ is bounded on $\R^N$,  for $\eps$ small enough, we have
$$
\uep(x) \geq \phi(x) +\frac 3 4 \gamma_r  \geq \phi(x) + \eps \upcor( \frac{x}{\eps})  +\frac{\gamma_r}{2} \quad\hbox{for all   }x \in  \partial B(\bar{x}, r) \; .
  $$ 
 
Choosing now $\gamma, r$ small enough, we can have, at the same time, (i) $\phi(\cdot)+ \eps \upcor( \frac{\cdot}{\eps})+\gamma$ is a subsolution of    \eqref{Bellman-OmH+} in $B(\bar{x},r)$ and (ii) $\gamma\leq\gamma_r /2$ in order that this function is less than $\uep$ on the boundary of this ball.

Therefore, applying the comparison result  (cf. Theorem~\ref{Umeno-comparisonPalle} (v)), we obtain
$$
\uep(x)   \geq \phi(x) + \eps \upcor( \frac{x}{\eps})  + \frac{\gamma}{2} \quad  \hbox{for all  } x \in   B(\bar{x}, \bar{r})  \; ,
  $$  
and letting $\eps \ds 0$  we have 
$$
\up(x) \geq  \phi(x) + \frac{\gamma}{2} \quad  \hbox{for all  } x \in   B(\bar{x}, \bar{r})  \; .
$$
Taking $x =\bar{x}$ we get $\up(\bar{x}) \geq \phi(\bar{x}) + \frac{\gamma}{2}$ which is in contradiction with \eqref{ipminstrettoSuper-up} and the proof is completed. 

\bigskip
 
\noindent {\bf Step 2  : $\up$  is a subsolution of \eqref{sys:limiteUpiu} --- the case when $b$ and $l$ do not depend on the first variable.}

In this step, we write for simplicity $b(y,a)$, $l(y,a)$ and $\mathcal{B}(y)$ since there is not dependence on the first variable. Let    $\phi$ be a $C^1$ function and   $\bar{x}$ a local (strict) maximum point for $\up -\phi$ such that $\up(\bar{x})-\phi(\bar{x})=0$, i.e
\begin{equation}  \label{max-upiu}
\exists \bar{r} >0  \:  \mbox{ such  that  } \up(x)-\phi(x) <0   \quad  \hbox{for all  } x \in {B(\bar{x},\bar{r})}   \mbox{ and } \up(\bar{x})-\phi(\bar{x})=0.
 \end{equation}
Our aim is to prove that
  \be  \label{tesi:subsol-upiu}
  \lambda  \phi(\bar{x}) +\bar{H}^+(\bar{x},D \phi(\bar{x}))  \leq 0\,,
 \ee 
 where $\Hbarp(\bar{x}, D\phi(\bar{x}))$  is defined  as in  Theorem   \ref{cell-Upiu-esistenza}. 
We argue by contradiction, assuming that 
   \be  \label{tesi:subsol-upiuASSURDA}
\exists \eta >0  \:  \mbox{ such  that   }   \lambda  \phi(\bar{x}) +\bar{H}^+(\bar{x},D \phi(\bar{x}))  \geq   \eta >0.
 \ee

Let $\upcor(y)$ be the   subsolution  related to  $\Hbarp(\bar{x}, D\phi(\bar{x}))$ fulfilling 
 \eqref{H+v}, therefore we have for any $y\in \R^N$ and $\tau\geq 0$ 
\begin{equation} \label{vdixbarrra} 
\upcor(y)=
 \inf_{(\trajYy,\a) \in \Taxbarrareg}  \left\{ \int_0^\tau \Big( l\big(\trajYy(t),a(t) \big) +b(\trajYy(t),a(t) ) \cdot D\phi(\bar{x})
    +\bar{H}^+(\bar{x},D\phi(\bar{x}))\Big)dt+ \upcor(\trajYy(\tau)) .
    \right\}
\end{equation}
Let $\eps>0$. In order to write explicitly $\upcor(\frac{\xb}{\e})$ we consider, for $y=\xb/\eps$,
any regular trajectory $\trajYy$ satisfying 
$\dot \trajYy (t)\in\mathcal{B}\big(\trajYy (t)\big)$, with $\trajYy(0)=\xb/\eps$.
Setting $X_\eps(t):=\eps \trajYy (\frac{t}{\eps})$, we obtain a solution of the differential inclusion
  \be  \label{dinupX}
  \dot{X}_\eps(t)\in \mathcal{B} \apt \frac{ X_\eps(t)}{\eps} \cht\,, \quad X_\eps(0)=\xb\,.
 \ee
We rewrite \eqref{vdixbarrra} by using the trajectories 
$\big(\trajYy(t),a(t)\big)$ under the form $\big(X_\eps(\eps t)/\eps,a(t/\eps)\big)$.
After a scaling in time, and rewriting $a(t)$ any control of the form $a(t/\eps)$, we arrive at
\begin{equation}   \label{convup3.0}
\begin{aligned}\eps  \upcor(\frac{\xb}{\eps})= 
\inf_{ (X_\e , \a )   \in \Taxbarrareg} \Bigg\{   \int_0^{ \e \tau} 
 \Big( l\big(\frac{X_\eps(t)}{\eps},\a(t)\big) +b\big(\frac{X_\eps(t)}{\e},\a(s)\big)
  \cdot D \phi(\bar{x}) 
 &  + \bar{H}^+(\bar{x}, D \phi(\bar{x}))  \Big)  dt \\
 & + \eps \upcor(\frac{ X_\eps  (\tau \e)}{\e}) \Bigg\}\,.
\end{aligned}
 \end{equation}
Our aim is now to prove that the function $w_\eps(x):=\phi(x)+\eps\upcor(\frac{x}{\eps})$ {\it{almost fulfills a super-optimality principle}} for $x=\bar{x}$ more precisely :

\begin{lem}\label{keypointconvpPIU} 
For each $\eps >0$ and $\bar t >0$   we have
\begin{equation} \label{DPPsuperupiu} 
w_\eps(\xb)  \geq  \inf_{(X_\e , \a )   \in \Taxbarrareg}  
\apg     
\int_0^{ \bar t}  l(\frac{X_\eps(t)}{\eps},\a(t)) e^{-\lambda t}  dt + w_\eps(X_\eps(\bar t))e^{-\lambda \bar t}    \chg + \eta \bar t+\eps O(\bar t) + o(\bar t)
\end{equation}
where $\eta$ is given by \eqref{tesi:subsol-upiuASSURDA} and the $o(\bar t)$ is uniform with respect to $\eps>0$.
\end{lem} 

{\bf Proof}. Since $\phi \in C^1(\R^\N)$ we have for each $\eps >0$ and $\bar t >0$   
 \be
 \phi(\xb)= \phi({X}_\e( \bar t )  ) e^{-\lambda \bar t } -
  \int_0^{\bar t} D \phi({X}_\e(t)) \cdot b(\frac{{X}_\eps(t)}{\e},\a(t)) e^{-\lambda t} dt  +
  \lambda  \int_0^{\bar t}   e^{-\lambda t} \phi( {X}_\e(t)) dt 
 \ee 
choosing $\tau=\frac{\bar t }{\eps}$ in \eqref{convup3.0}
  \begin{align*} 
 w_\eps(\xb)= & \inf_{ (X_\e , \a )   \in \Taxbarrareg}  \Bigg\{  \int_0^{ \bar t} 
 \apt l(\frac{X_\eps(t)}{\eps},\a(t)) +b(\frac{X_\eps(t)}{\e},\a(t)) \cdot D \phi(\bar{x})   + \bar{H}^+(\bar{x}, D \phi(\bar{x}))  \cht  dt + \eps \upcor(\frac{ X_\eps  (\bar t)}{\e})\\
 &   \qquad + \ap  \phi({X}_\e(\bar t)) e^{-\lambda \bar t} -
  \int_0^{\bar t} D \phi({X}_\e(t)) \cdot b(\frac{{X}_\eps(t)}{\e},\a(t)) e^{-\lambda t}  dt  +  \lambda  \int_0^{\bar t}   e^{-\lambda t} \phi( {X}_\e(t)) dt  \chg  \\
   = & \inf_{ (X_\e , \a )   \in \Taxbarrareg}  \Bigg\{   \int_0^{ \bar t} 
       l(\frac{X_\eps(t)}{\eps},\a(t))  e^{-\lambda t} dt +  w_\eps({X}_\e(\bar t)  )
e^{-\lambda \bar t} + \int_0^{\bar t}  \big(\bar{H}^+(\bar{x}, D \phi(\bar{x})) + \lambda \phi(\xb)  \big) dt\\ 
& \qquad +  \eps V^+(\frac{ X_\eps  (\bar t)}{\e})(1-e^{-\lambda \bar t})
+ \int_0^{ \bar t} (1-e^{-\lambda  t}) l(\frac{X_\eps(t)}{\eps},\a(t))  dt \\ & \qquad +
\int_0^{ \bar t}   b(\frac{X_\eps(t)}{\e},\a(t)) \cdot D \phi(\bar{x}) -   D \phi({X}_\e(t)) \cdot b(\frac{{X}_\eps(t)}{\e},\a(t))  e^{-\lambda t}  dt\Bigg\}\,.
\end{align*}
Notice that since $V^+$ is bounded, $|\eps V^+(X_\eps(\bar t)/\e)(1-e^{-\lambda \bar t})|\leq C\eps\bar t=\eps O(\bar t)$ (note that the $O(\bar t)$ is  uniform in $\eps$). In order to estimate the two last terms, we use the fact that if $g$ is bounded, then  
$$\bigg|\int_0^{\bar t} g(t) dt - \int_0^{\bar t} g(t) e^{-\lambda t} dt\bigg|\leq
 \|g\|_\infty(\bar t)^2 = o(\bar t)\,,$$
where the $o(\bar t)$ only depends on $\|g\|_\infty$. Hence, since the trajectory $X_\eps$ is continuous,  $b$ and $l$ are bounded and $\phi \in C^1(\R^\N)$, we have 
\begin{align*} w_\eps(\xb) = &
\inf_{ (X_\e , \a )   \in \Taxbarrareg}  \Bigg\{   \int_0^{ \bar t} 
    l(\frac{X_\eps(t)}{\eps},\a(t))  e^{-\lambda t} dt +  w_\eps({X}_\e(\bar t)) e^{-\lambda \bar t}  + \int_0^{\bar t}  \big(\bar{H}^+(\bar{x}, D \phi(\bar{x})) + \lambda \phi(\xb)  \big) dt \Bigg\} + o(\bar t)\\
\geq & \inf_{ (X_\e , \a )   \in \Taxbarrareg}  \Bigg\{   \int_0^{ \bar t} 
    l(\frac{X_\eps(t)}{\eps},\a(t))  e^{-\lambda t} dt +  w_\eps({X}_\e(\bar t)) e^{-\lambda \bar t}  \Bigg\} + \eta \bar t+\eps O(\bar t) + o(\bar t)\,,
\end{align*}
which gives the result.
\hfill $\Box$

\bigskip

We consider now  the DPP for the function $\uep$ at point $ \xb $ and time $\bar t$  
\begin{equation}  \label{DPPUeps}
   \uep (\xb)=   \inf_{ (X_\e , \a )   \in \Taxbarrareg} 
    \left\{ \int_0^{ \bar t}   l(\frac{ X_\eps(t) }{\e},a(t) \big)  e^{-\lambda t} dt
    +  \uep (X_\eps( \bar t))  e^{-\lambda  \bar t} \right\}
\end{equation}
and combine it with \eqref{DPPsuperupiu} to get
 \begin{equation*} 
\uep (\xb)- w_\eps(\xb) \leq  \sup_{ (X_\e , \a )   \in \Taxbarrareg}  
 \Bigg\{  \Big( \uep (X_\eps( \bar t))
  - w_\eps( X_\eps( \bar t)) \Big ) e^{-\lambda\bar t} \Bigg\} - \eta \bar t 
  +\eps O(\bar t) + o(\bar t).  
\end{equation*} 
Therefore, using again that $V^+$ is bounded (and $\bar t$ can be chosen, say, less than 1), we get
\begin{equation} \label{eqxbarre}
\uep (\xb)- \phi(\xb) \leq  \sup_{ (X_\e , \a )   \in \Taxbarrareg}  
 \Bigg\{ \Big( \uep (X_\eps( \bar t))
  - \phi( X_\eps( \bar t)) \Big ) e^{-\lambda\bar t} \Bigg\} - \eta \bar t +\eps O(\bar t)+o(\bar t)+o_\eps(1).  
\end{equation} 
We choose now $\bar t \leq r/(2\|b\|_\infty)$ in order that the trajectory $X_\eps(t)$ belongs to $B(\xb, \bar r)$ for all $t\in[0,\bar t]$. In particular, since $X_\eps(\bar t)\in B(\xb, \bar r)$, \eqref{max-upiu} implies that
$$ \uep (X_\eps( \bar t))-\phi( X_\eps( \bar t) )\leq o_\eps(1)\,.$$
We then pass to the limit as $\eps\to0$ in \eqref{eqxbarre} and get
\begin{equation*} 
0=\up (\xb)- \phi(\xb) \leq  - \eta \bar t +o(\bar t)\,,
\end{equation*} 
which is a contradiction for $\bar t$ small enough. Hence the proof is complete.

\bigskip
 
\noindent {\bf Step 3  : $\up$  is a subsolution of \eqref{sys:limiteUpiu} --- the general case.}

In the above proof, the key fact was that the change we made on the trajectories $\trajYy$ lead us to \eqref{dinupX} which is exactly the dynamic for the control problem which gives $\uep$. On the contrary, when $b$ depends on the first variable, this change is going to provide a dynamic $b(\bar x, \frac{\cdot}{\eps}, a(t))$ instead of $b(\cdot, \frac{\cdot}{\eps}, a(t))$: if $b$ was continuous, we could handle this difference but here it may change the times when $\H$ is reached or left by the trajectories and we cannot compare the control problem for $\uep$ and $w_\eps$.

In order to overcome this difficulty, we introduce 
a $\kappa$-sequence of  problems for $\kappa\in]0,1]$ where, for $\kappa$ fixed, the dynamics are constant with respect to the slow variable and the new trajectories can be interpreted as 
a subset  of ${\cal T}^{\rm reg}$. For $\kappa$ fixed we shall use Step
2 which provides a modified corrector $V^+_\kappa$, solution of a suitable ergodic problem and a modified ergodic constant $\bar{H}^+_\kappa$.
Then we conclude letting $\kappa \rightarrow 0$, indeed, thanks to the stability properties proved in \cite{BBC2} the ergodic constants are stable.\\
As in Step 2 let $\phi$ be a $C^1$ function and   $\bar{x}$ a local (strict) maximum point for $\up -\phi$ such that $\up(\bar{x})-\phi(\bar{x})=0$, i.e
\begin{equation}  \label{max-upiuu}
\exists \bar{r} >0  \:  \mbox{ such  that  } \up(x)-\phi(x) <0   \quad  \hbox{for all  } x \in {B(\bar{x},\bar{r})}   \mbox{ and } \up(\bar{x})-\phi(\bar{x})=0.
 \end{equation}
Our aim is to prove (\ref{tesi:subsol-upiu}), i.e.
\[
  \lambda  \phi(\bar{x}) +\bar{H}^+(\bar{x},D \phi(\bar{x}))  \leq 0 .
 \]
In order to define the approximate corrector, we first introduce, for $0<\kappa \ll 1$ and $i=1,2$, the sets
$$\mathbf{BL}_{i\kappa}(\xb,y):=\bigcap_{|z-\xb|\leq\kappa}\mathbf{BL}_i(z,y)\,,$$
where
$$\mathbf{BL}_i(z,y):= \big\{ (b_i(z,y,\alpha_i),l_i(z,y,\alpha_i)) : \alpha_i \in
    A_i \big\}\,,\quad \text{for any } (z,y)\in\R^N\times\Omega_i\,.  $$

With such definitions, we can build $\mBL_\kappa$ as we built $\mB$ for \eqref{def:traj} and we are interested in the differential inclusion 
\begin{equation}\label{dif.inclusion.kappa}
    \frac{d}{dt}\big(\trajYy (t),\trajLy (t)\big)\in\mBL_\kappa(\xb,\trajYy(t))\,,\text{ with }(\trajYy (0),\trajLy (0))=(y,0)\,.
\end{equation}
Notice that solutions of this differential inclusion are couples $\big(Y_y,L_y\big)$, but that
$\mathcal{BL}_\kappa$ does not depend on $L_y$. 

So, despite we are not exactly 
in the framework of Theorem~\ref{def:dynAssume},
there is no difficulty to solve differential inclusion
\eqref{dif.inclusion.kappa}, using [H0],
[H1], [H2], since the set-valued map $\mBL_\kappa$ is upper
semi-continuous with convex compact images. And, as in
Theorem~\ref{def:dynAssume}, for each solution  $(\trajYy (\cdot),\trajLy
(\cdot))$ of \eqref{dif.inclusion.kappa} and for each function $e\in
L^\infty((0,+\infty))$ with $\| e(t)\| \leq 1$ a.e.,  there exists a control
$\a(\cdot)=\big(\alpha_1(\cdot),\alpha_2(\cdot),\mu(\cdot)\big) \in \Ac$ such that
\begin{equation}\label{eq:repY}
\begin{aligned}
\dot\trajYy (t)  =& b_1\big(\xb +\kappa e (t), \trajYy (t),\alpha_1(t)\big)\mathds{1}_{\apg \trajYy (t) (t) \in  \Omega_1 \chg }+
b_2\big(\xb +\kappa e (t), \trajYy (t),\alpha_2(t)\big)\mathds{1}_{\apg  \trajYy(t) \in \Omega_2 \chg } \\[2mm]
&+ b_\H\big(\xb +\kappa e (t), \trajYy (t),\a(t)\big)\mathds{1}_{\apg  \trajYy (t) \in  \H \chg }\,, 
\end{aligned}
\end{equation}
\begin{equation}\label{eq:repL}
\begin{aligned}
\dot\trajLy (t) =& l_1\big(\xb +\kappa e (t), \trajYy (t),\alpha_1(t)\big)\mathds{1}_{\apg \trajYy (t) (t) \in  \Omega_1 \chg }+
l_2\big(\xb +\kappa e (t), \trajYy (t),\alpha_2(t)\big)\mathds{1}_{\apg  \trajYy(t) \in \Omega_2 \chg } \\[2mm]
&+ l_\H\big(\xb +\kappa e (t), \trajYy (t),\a(t)\big)\mathds{1}_{\apg  \trajYy (t) \in  \H \chg }\,.
\end{aligned}
\end{equation}
A key remark here is that, for the associated control problem, the running cost is going to be given by $L_y(\cdot)$ and, if we fix the solution of the differential inclusion, is {\em independent of the choice of the control} $\a(\cdot)$. However we have to define the set of regular  trajectories  $(\Tayarrareg)_\kappa$ as for the original problem and we say that $(\trajYy (\cdot),\trajLy (\cdot))$ is a regular trajectory if there exists function $e\in L^\infty((0,+\infty))$ with $|e(t)|\leq 1$ a.e. such that $(\trajYy (\cdot),\a(\cdot))$ is a regular trajectory. 

By analogy with \eqref{elletilde3}, we replace 
$$
\tilde{l}\big(x,p,\trajYy (t),a(t) \big)\quad\text{by}\quad  \dot L_y (t) + \dot Y_y (t)\cdot p 
$$
in \eqref{def-costo-y}, for any $\big(Y_y(\cdot),L_y(\cdot)\big)$
satisfying the dynamics \eqref{eq:repY}--\eqref{eq:repL} for some
controls $e\in L^\infty(0,\infty)$, $\|e(\cdot)\|\leq1$ and $a\in
\mathcal{A}$\,.

Following the $\rho$-problem construction as in the proof of
Theorem~\ref{cell-Upiu-esistenza}, we obtain a
Lipschitz continuous periodic function $V^+_\kappa$ (the approximate
corrector) such that, for any $\tau >0$
$$
V^+_\kappa (y)=
 \inf_{(\trajYy,\trajLy) \in (\Tayarrareg)_\kappa}  \left\{ \trajLy(\tau)+ (\trajYy(\tau)-y) \cdot D\phi(\bar{x})
    +\tau \bar{H}^+_\kappa (\bar{x},D\phi(\bar{x}))+ V^+_\kappa (\trajYy(\tau))
    \right\}\; .
$$
To do so, it is worth pointing out that,  [H2] being uniform with respect to $x$ and $y$, the trajectory $\trajYy$ is still controllable.

%

In order to prove (\ref{tesi:subsol-upiu}), we are going to show that
\begin{equation}\label{eq:ergo.kappa}
    \lambda\phi(\xb)+\bar{H}^+_\kappa(\xb, D\phi(\xb))\leq 0\, ,
\end{equation}
and then we will prove that $\bar{H}^+_\kappa(\xb, D\phi(\xb)) \to \bar{H}^+(\xb, D\phi(\xb))$.

We argue by contradiction assuming that $\lambda\phi(\xb)+\bar{H}^+_\kappa(\xb, D\phi(\xb))\geq \eta$ for some $\eta >0$.

We introduce the function $w_{\eps,\kappa}$ defined by
$w_{\eps,\kappa} (x):=\phi(x)+\eps V^+_\kappa (\frac{x}{\eps})$ and we first examine the properties of $\eps V^+_\kappa (\frac{x}{\eps})$.

For any solution $(Y_y,L_y)$ of \eqref{dif.inclusion.kappa}, if we set $Z_\eps(t)= (Z^1_\eps(t),Z^2_\eps(t)):=(\eps \trajYy (t/\eps),\eps \trajLy (t/\eps))$ with $y = \xb /\eps$, this trajectory solves
$$
 \dot Z_\eps(t) \in\mBL_\kappa\Big(\xb,\frac{Z^1_\eps(t)}{\eps}\Big)\,,
 \text{ with }\big(Z^1_\eps(t) ,Z^2_\eps(t)\big) (0)=(\xb,0)\,.
 $$
Our hypotheses on the dynamics imply that for any $\kappa$, there exists
$\bar t=\bar t(\kappa)>0$ such that for all
$Z_\eps\in(\Taxbarrareg)_\kappa$ and all $t\in[0,\bar t]$,
$|Z^1_\eps(t)-\xb|<\kappa$. By the definition of $\mBL_\kappa$ (choosing $e(t)=\frac{Z^1_\eps(t)-\bar x}{\kappa}$ ), this
implies that, for such~$t$
$$ \mBL_\kappa\Big(\xb,\frac{Z^1_\eps(t)}{\eps}\Big) 
\subset \mBL\Big(Z^1_\eps(t),\frac{Z^1_\eps(t)}{\eps}\Big) \; .$$
This key inclusion property means that the $Z^1_\e$-trajectories can be
seen as particular regular $X_\eps$-trajectories on $[0,\bar t]$ of 
$\dot X_\eps(t)\in\mathcal{B}\big(X_\eps,\frac{X_\eps}{\eps}\big)$.
 Hence, using \eqref{eq:repY}--\eqref{eq:repL} in the definition of $V^+_\kappa$
and taking the infimum on a bigger set
we conclude that, for any $t\in [0, \bar t]$
$$
\eps V^+_\kappa (\frac{\xb}{\eps}) \geq  \inf_{(X_\e , \a )   \in \Taxbarrareg}  \left\{ \int_0^{t} \left(
	\tilde{l} \big(X_\eps(s),D\phi(\bar{x}), \frac{X_\eps(s)}{\eps},\a(s)\big)+\bar{H}^+_\kappa (\bar{x},D\phi(\bar{x}))\right)ds+ \eps V^+_\kappa (\frac{X_\eps(t)}{\eps})
    \right\}\; .$$

In the previous formula, since the integrand is bounded, we also have
$$
\begin{aligned}
\eps V^+_\kappa (\frac{\xb}{\eps}) \geq  \inf_{(X_\e , \a )   \in \Taxbarrareg} & \bigg\{ \int_0^{t} \left(
\tilde{l} \big(X_\eps(s),D\phi(\bar{x}), \frac{X_\eps(s)}{\eps},\a(s)\big)+\bar{H}^+_\kappa (\bar{x},D\phi(\bar{x}))\right)e^{-\lambda s}ds\\
&+ \eps V^+_\kappa (\frac{X_\eps(t)}{\eps})e^{-\lambda t}\biggr\} + o(t)+ \eps O(t)\; .
    \end{aligned}$$
Finally, considering the function $s\mapsto \phi(X_\eps(s))e^{-\lambda s}$, we have
\begin{eqnarray}
\phi(X_\eps(t))e^{-\lambda t} - \phi(\xb) &=& \int_0^t (D\phi(X_\eps(s))\cdot \dot X_\eps(s)-\lambda \phi(X_\eps(s)))e^{-\lambda s}ds\\
&=&  \int_0^t \left(D\phi(X_\eps(s))\cdot b(X_\eps(s),\frac{X_\eps(s)}{\eps},a(s))-\lambda \phi(X_\eps(s))\right)e^{-\lambda s}ds
\end{eqnarray}
 and, taking in account the facts that $\phi$ is smooth and $b$ is uniformly bounded, we have
$$
\int_0^t D\phi(\xb)\cdot b(X_\eps(s),\frac{X_\eps(s)}{\eps},a(s))e^{-\lambda s}ds = \phi(X_\eps(t))e^{-\lambda t} - \phi(\xb) + \int_0^t \lambda \phi(\xb)e^{-\lambda s}ds + o(t)\; .
$$
We deduce from these properties that, for any $\eps>0$ and $t \leq \bar t$
$$
\begin{aligned} 
w_{\eps,\kappa}(\xb)  \geq  \inf_{(X_\e , \a )   \in \Taxbarrareg}  
\biggl\{   
& \int_0^{t}\left( l\big(X_\eps(s),\frac{X_\eps(s)}{\eps},\a(s)\big) +\bar{H}^+_\kappa(\xb, D\phi(\xb))+ \lambda \phi(\xb)\right)e^{-\lambda s} 
ds \\&+ w_{\eps,\kappa}\big(X_\eps(t)\big)e^{-\lambda t}    \biggr\} +o(t)+\eps O(t)\,.
\end{aligned}
$$
Finally we use that $\bar{H}^+_\kappa(\xb, D\phi(\xb))+ \lambda \phi(\xb) \geq \eta >0$
\begin{equation}\label{DPPsuperupiuGEN} 
w_{\eps,\kappa}(\xb)  \geq  \inf_{(X_\e , \a )   \in \Taxbarrareg}  
\biggl\{   
\int_0^{t}\left(  l\big(X_\eps(s),\frac{X_\eps(s)}{\eps},\a(s)\big)\right)e^{-\lambda s} 
ds \\ + w_{\eps,\kappa}\big(X_\eps(t)\big)e^{-\lambda t}    \biggr\} + \eta t +o(t)+\eps O(t)\,.
\end{equation}

In order to conclude, we consider  the DPP for the function $\uep$ at point $ \xb $ and time $t$ 
\begin{equation*} 
   \uep (\xb)=   \inf_{ (X_\e , \a )   \in \Taxbarrareg} 
    \left\{ \int_0^{t}   l(X_\eps(s),\frac{ X_\eps(s) }{\e},a(s) \big)  e^{-\lambda s} ds
    +  \uep (X_\eps(t))  e^{-\lambda t} \right\}\; .
\end{equation*}
Since combining with  \eqref{DPPsuperupiuGEN}  we obtain 
\begin{equation*}
\uep (\xb)- w_{\eps,\kappa}(\xb) \leq    \sup_{ (X_\e , \a )   \in \Taxbarrareg} 
  \bigg\{ U^+_\eps\big(X_\eps(t)\big)-w_{\eps,\kappa}\big(X_\eps( t)\big)
  \bigg\} - \eta t+\eps O( t) +o( t).   
\end{equation*} 
Arguing now exactly as in Step 2, we have the contradiction and we have proved that, for any $\kappa>0$, \eqref{eq:ergo.kappa} holds.\\

The final step consists in passing to the limit as $\kappa\to0$, thanks to the
\begin{lem}\label{lem:H+kappa}
    When $\kappa\to0$, $\bar{H}^+_\kappa(\xb,p)\to \bar H^+(\xb,p)$ for any fixed
    $\xb,p$.
\end{lem}

For the sake of clarity we prove  this lemma below and  we provide now  the conclusion of the proof.  
Passing  to the limit in \eqref{eq:ergo.kappa} we  get 
$\lambda\phi(\xb)+\bar H^+(\xb, D\phi(\xb))\leq0$. Hence, $\up$ is a subsolution of the ergodic problem  \eqref{sys:limiteUpiu}. Combined with Step 1, the proof of Theorem~\ref{Upiu-convergenza} is complete.
\hfill $\Box$

{\bf Proof of Lemma~\ref{lem:H+kappa}}.
Since the $\{V^+_\kappa(\cdot)\}_\kappa$ are equi-Lipschitz (by [H2]) and 
periodic, after normalizing by $V^+_\kappa(0)=0$ they are 
uniformly bounded. Applying Ascoli-Arzela Theorem, up to the extraction of a subsequence, we may assume that $V^+_\kappa\to V$ locally uniformly, for  some Lipschitz periodic continuous function $V$ (which, a priori, may not be $V^+$ even up to an additive constant). Similarly, we can assume 
that for any fixed $(\xb,p)$, $\bar{H}^+_\kappa(\xb,p)\to C^+(\xb,p)$ for some 
constant $C^+$.

Since, for any $\tau \geq 0$ we have
$$
V^+_\kappa (y)=
 \inf_{(\trajYy,\trajLy) \in (\Tayarrareg)_\kappa}  \left\{ \trajLy(\tau)+ (\trajYy(\tau)-y) \cdot p    +\tau \bar{H}^+_\kappa (\bar{x},p)+ V^+_\kappa (\trajYy(\tau))
    \right\}\,,
$$
for each $\kappa>0$ there exists a regular $\kappa$-trajectory $(\trajYy^\kappa,\trajLy^\kappa)$ such that 
this infimum is attained. To pass to the limit, it is clear that we have
a subsequence of trajectories which converges uniformly to a trajectory
$(Y,L)$ of the limit problem, thanks to $[H0]-[H1]$. However, we need the 
limit trajectory $Y$ to be regular. 

So, in order to use Lemma 5.3 in \cite{BBC2}, we first remark that, if
$\mBL$ is built in the same way as $\mBL_\kappa$, using $\mathbf{BL}_i$
instead of $\mathbf{BL}_{i\kappa}$, then any solution of the
$\mBL_\kappa$--differential inclusion is a solution of the
$\mBL$--differential inclusion since $\mBL_\kappa(\xb,y)\subset
\mBL(\xb,y)$.

Thus, if the subsequence $(\trajYy^{\kappa_n},\trajLy^{\kappa_n})$
converges to $(\trajYy^{*},\trajLy^{*})$, then
$(\trajYy^{*},\trajLy^{*})$ solves the $\mBL$--differential inclusion.
Moreover, since the trajectories $(\trajYy^{\kappa_n},\trajLy^{\kappa_n})$
are regular, there exists $(e^{\kappa_n},a^{\kappa_n})$ such that
$(\trajYy^{\kappa_n},\trajLy^{\kappa_n})$ can be interpreted as a
regular trajectory in $\R^N \times \R$, associated with a partition
$(\Omega_1\times \R) \cup (\Omega_2\times \R)\cup(\H\times \R)$,
$e^{\kappa_n}$ being interpreted as part of the control.

The assumptions of  \cite[Lemma 5.3]{BBC2} are then fulfilled and we deduce
that $(\trajYy^{*},\trajLy^{*})$  is a regular trajectory. This means
that there exists a control $a^*$ (of course, there is no $e$ playing a
role at the limit) such that the limit trajectory is regular:
$(\trajYy^*,a^*)\in \Tayarrareg$, and we
have
$$V(y)= 
  \int_0^\tau \Big( l\big(\bar{x},\trajYy^*(t),a^*(t) \big) 
 +b(\xb,\trajYy^*(t),a^*(t) ) \cdot p+ 
 C^+(\xb,p)\Big)dt + V(\trajYy^*(\tau))\,.
$$
Hence for any $\tau>0$
$$V(y)\geq
 \inf_{(\trajYy,\a) \in \Tayarrareg}  \left\{ \int_0^\tau \Big( l\big(\bar{x},\trajYy (t),a(t) \big) +b(\bar{x},\trajYy(t),a(t) ) \cdot p    +C^+(\bar{x},p)\Big)dt+ V(\trajYy(\tau))
    \right\}\,.
$$

Now we prove the converse inequality. To do so we notice that the function 
$V^+_\kappa(y)-H^+_\kappa(\xb,p)\cdot t$ 
is a subsolution of the evolution problem  
$$w_t+\mathbf{H}^+_\kappa(\bar x,y,Dw+p)=0\;\hbox{in  }\R^N \times (0,+\infty),  \quad w(y,0)=V^+_\kappa(y)\;\hbox{in  }\R^N \,.$$
where $\mathbf{H}^+_\kappa$ is defined as $\mathbf{H}^+$, but with the 
$\mathbf{B}_{i\kappa}$.  
Arguing as above 
we can pass to the limit as $\kappa\to0$ using  the
stability result for $\mathbf{H}^+$  \cite[Theorem 5.1]{BBC2} and  we obtain
that $V(y)-C^+(\xb,p)\cdot t$ is a subsolution of 
$$w_t+\mathbf{H}^+(\bar x,y,Dw+p)=0\;\hbox{in  }\R^N \times (0,+\infty),  \quad w(y,0)=V(y)\;\hbox{in  }\R^N \,.$$
Hence, writing down the maximal (sub-)solution of this problem, we 
obtain that for any $y$ and $\tau$,
$$V(y)-C^+(\bar{x},D\phi(\bar{x}))\cdot  \tau \leq 
\inf_{(Y_y,\a) \in \mathcal{T}^{\rm reg}_y}  \left\{ \int_0^\tau \Big( 
    l\big(\bar{x},Y_y(t),a(t) \big) +b(\bar{x},Y_y(t),a(t) ) \cdot p
    \Big)dt+ V(Y_y(\tau))
    \right\}\,.
$$
Combining both inequalities for $V(y)$, we get equality and the characterization of $\bar H^+(\xb,p)$ in Theorem~\ref{cell-Upiu-esistenza} implies that
necessarily, $C^+(\bar{x},p) =\bar H^+(\xb, p)$.
\hfill $\Box$

\begin{rem}
Actually, in Step 3 we only need the inequality $C^+(\bar{x},p) \geq\bar H^+(\xb, p)$ to conclude that 
$U^+$ is a subsolution. However, knowing that $H^+_\kappa(\bar{x},p) \to \bar H^+(\xb, p)$ is of an independent interest.
\end{rem}


\section{The 1-D case: an example}\label{sect:oneD}

In this section we present an $1$-d example showing that $\Hbarm$ and $\Hbarp$ can be indeed different. We refer to 
Lions, Papanicolaou and Varadhan \cite{LPV}, Concordel \cite{Conc1, Conc2} and Namah and Roquejoffre \cite{NR} for explicit computations of effective Hamiltonians : despite our aim is not really to compute explicitly $\bar{H}^\pm (x,p)$, our arguments are inspired by these works.

We choose $\displaystyle \Omega_1=\bigcup_{k\in\Z} ]0,1[+2k$ and $\displaystyle \Omega_2=\bigcup_{k\in\Z}]1,2[+2k$. In this context, $\umcor$, $\upcor$ are $2$-periodic functions and from (\ref{H+v}) (and its analogue for $\umcor$), we have
\begin{equation}\label{DPP-pm}
V^\pm(y_0)=
 \inf_{(\trajyoyo,\a)}  \left\{\int_0^t   l\big(x,\trajyoyo(s),a(s) \big)  ds
   +(\trajyoyo(t)-y_0)  \cdot p
    +V^\pm(\trajyoyo(t))
    \right\}+\bar{H}^\pm (x,p) t \; ,
\end{equation}
where we have used that $\trajyoyo$ solves (\ref{def:trajyo}). Of course, the admissible controls are different for $\umcor$, $\upcor$.

Since we are interested in the ergodic constant, we can choose $y_0=0$ and we also consider large $t$. For the optimal trajectories we have two possible cases.

\underline{\it Case 1: } for some $t>0$, there exists $0<\bar{t} \leq t$ such that $Y_0(\bar{t})=\pm 2$. 

We point out that, by the Dynamic Programming Principle, this case is equivalent to $|Y_0(s)| \geq 2$ for some $s \in (0,t]$ and clearly this happens, for large $t$, if $|p|$ is large since the term $Y_0(t) \cdot p$ is playing a more important role in the minimization process than the $\int_0^t   l\big(x,Y_0(s),a(s) \big)  ds$--one.

Recalling the periodicity of $V^\pm$, $V^\pm(\pm2)=V^\pm(0)$ and we obtain from (\ref{DPP-pm})
$$ 
\bar{H}^\pm (x,p) = - \frac{1}{\bar t} \int_0^{\bar t}   l\big(x,Y_0(s),a(s) \big)  ds
   \pm \frac{2}{\bar t}  p  = - \frac{1}{\bar t} \int_0^{\bar t}   l\big(x,Y_0(s),a(s) \big)  ds
   + \frac{2}{\bar t}  |p|\; .
$$
This gives the behavior of $\bar{H}^\pm (x,p)$ and one can be more precise about the different terms of this equality since $\bar t, Y_0, a$ solves the control problem
$$ \inf_{(Y_0,\a,t)}  \left\{\int_0^t   l\big(x,Y_0(s),a(s) \big)  ds\right\}\; ,$$
where the infimum is taken on all the trajectories such that $|Y_0(t)|=2$.

\underline{\it Case 2: }  For any $t$ and $s\in (0,t)$,  $|Y_0(s)| \leq 2$.

In that case, since $V^\pm$ are bounded and the optimal trajectory too, we divide (\ref{DPP-pm}) by $t$ and letting $t$ tend to $+\infty$ we obtain
 \begin{equation}  \label{formHp}
 \bar H^\pm (x,p)= \lim_{t \rightarrow +\infty}  \apt -  \inf_{(Y_0,\a)}  \big\{  \frac{1}{t}     \int_0^t   l\big(x,Y_0(s),a(s) \big)  ds
 \big\} \cht \; .
\end{equation}
Again we insist on the fact that the infimum is taken on different sets of controls for $\bar H^+$ and $\bar H^-$.

At this point, we consider the following example. 
Let $b_1(x,y,\alpha_1)=\alpha_1$ and $b_2(x,y,\alpha_2)=\alpha_2$, $A_1=A_2=[-1,+1]$ and
\[
l_1(x,y,\alpha_1)=\vert\alpha_1-\cos(\pi y)\vert +1-\vert\cos(\pi y)\vert
\]
\[
l_2(x,y,\alpha_2)=\vert\alpha_2+\cos(\pi y)\vert +1-\vert\cos(\pi y)\vert .
\]
If we consider the case $p=0$ then

\def\Tayoreg{{\cal T}^{\rm reg}_{y_0}}

 \[ \Hbarp(x,0)= \lim_{t \rightarrow +\infty}  \apt -  \inf_{(Y_0,\a) \in {\cal T}^{\rm reg}_{0}}  \big\{  \frac{1}{t}     \int_0^t   l\big(x,Y_0(s),a(s) \big)  ds  \big\} \cht 
\]
\[ \Hbarm(x,0)= \lim_{t \rightarrow +\infty}  \apt -  \inf_{(Y_0,\a) \in {\cal T}_{0}}  \big\{  \frac{1}{t}     \int_0^t  l\big(x,Y_0(s),a(s) \big)  ds  \big\} \cht.
\]

In the $\Hbarm(x,0)$-case, it is clear that the best strategy is to choose $\alpha_1=1$ and $\alpha_2=-1$. Indeed, while $l_1, l_2 \geq 0$, this strategy allows to stay at $0$ with a zero cost because $\bH (x,0,(1, -1, \frac 1 2))=0$ and $\lH  (x,0,(1, -1, \frac 1 2)):=\frac 1 2 l_1 (x,0,1) + \frac 1
2l_2(x,0,-1)=0$. Hence $ \Hbarm(x,0)= 0$.

But this trajectory is a ``singular'' one and cannot be used in the $\Hbarp(x,0)$-case. If only the ``regular'' strategies are allowed, then either the best strategy is to have a trajectory $Y_0$ which stays in $0$ but the controls have to satisfy
$\alpha_1\leq 0$ and $\alpha_2\geq 0$ so $\lH(x, 0, a)\geq1$ and this would give $\Hbarp(x,0)=
-1$ because $l_1(x,0,0)=1$,  $l_2(x,0,0)=1$. Or the best regular strategy consists in leaving $0$ but in this case, we  know from the Dynamic Programming Principle that the optimal trajectory has to be monotone (see \cite[Section 5.1]{BBC}, Lemma 5.4 in particular) and we can conclude that $\Hbarp(x,0) <0$. 

Indeed, either we are in a situation which is analogue to Case 1 above but we remark that, for the trajectory, the cost to go from $0$ to $1$ and then from $1$ to $2$ [or from $0$ to $-1$ and then from $-1$ to $-2$] is strictly positive. Or this trajectory remains bounded and therefore converges to some point in $[-2,2]$ but in this case it is easy to prove that the cost is larger than $1$ using that $\cos(\pi Y(s))$ is almost constant for large $s$ and that it is not possible to stay in $y=1$ with a "regular" trajectory.
 
Therefore we have shown that
\[
    \Hbarp(x,0)\not=\Hbarm(x,0)\,.
\]


%


\section{Appendix-Regularity results for the tangential Hamiltonian}  

In this section we can use some weaker hypothesis on $b$, we can replace [H0] by the following hypothesis. 

\begin{itemize}
\item[{[H0bis]}] 
  For $i=1,2$, $A_i$ is a compact metric space and $b_i : \R^N \times \R^N \times A_i  \ds \R^N$ is a continuous bounded function. 
More precisely, there exists  
$M_b >0$, such that for any $x \in \R^N$, $y \in \R^N$ and $\alpha_i \in A_i$, $i=1,2$,
$$ |b_i(x,y,\alpha_i) | \leq M_b \; . $$ 
 There exists a modulus $\omega_b(\cdot)$ such that, for any $x,z \in \R^N$, $y \in \R^N$ and $\alpha_i \in A_i$
$$ |b_i(x,y,\alpha_i)-b_i(z,y,\alpha_i)|  \leq \omega_{b}( x-z)\;\; , $$ 
and there exists $\bar{L}_{i}\in \R$ such that, for any $x\in \R^N$, $y,w \in \R^N$ and $\alpha_i \in A_i$
$$ |b_{i}(x,y,\alpha_i)-b_{i}(x,w,\alpha_i)|\leq \bar{L}_{i} |y-w|\;.$$ 
\end{itemize}

In order to obtain some regularity result for the tangential Hamiltonian we prove  the following useful property
\begin{lem}   \label{lemmacont}
Assume [H0bis], [H1], 
and [H2], let $\H$ be a $W^{2,\infty}$-hypersurface. 
Fix $(x,y), (z,w) \in \R^\N \times \H$. For each  control $a \in A_0(x,y)$, there exists a control $\tilde{a} \in A_0(z,w)$ such that 
\begin{equation}  \label{lemmaconteq}
|( \bH(x,y,a),\lH(x,y,a)) -  ( \bH(z,w,\tilde{a}),\lH(z,w,\tilde{a}))  | \leq C |  \omega_{b}( x-z)+ \omega_{l}( x-z) +| y-w|   | 
\end{equation}
where $C$ is an explicit constant depending on  $\bar{L}_i$, $\bar{L}_{i,l}$, $M_b$, $M_l$, $\delta$ introduced in [H0], [H1], [H2] and on the Lipschitz constant $L_n$  of the normal vector $ \nor_1$.
\end{lem}  
{\bf Proof.}
First we remark that if $a\in A$
\begin{equation} \label{lipafix}
|( \bH(x,y,a),\lH(x,y,a)) -  ( \bH(z,w,a),\lH(z,w,a))  | \leq   |  \omega_{b}( x-z)+ \omega_{l}( x-z) +| y-w|   \max_{i=1,2}(\bar{L}_i,\bar{L}_{i,l}) |  
\end{equation}
Let us consider a control $a \in \A_0(x,y)$, i.e. $\bH(x,y,a) \cdot \nor_1(y)=0$.\\We have two possibilities.  
First, if $\bH(z,w,a) \cdot \nor_1(w) =0$ the conclusion follows  easily using \eqref {lipafix}  because $a \in \A_0(z,w)$. \\
In the second case,  let us suppose that  $\bH(z,w,a) \cdot \nor_1(w) > 0$. (For the other sign the same argument will apply so we will not detail it.)
By the controllability assumption in [H2] there exists a control $a_1 \in A$ such that  $ \bH(z,w,a_1) \cdot \nor_1(w)=- \delta$ . 
We set now 
$$
 \bar{\mu}:= \frac{ \delta  }{ \bH(z,w,a) \cdot \nor_1(w)  +  \delta  }, 
$$
since $ \bar{\mu} \in ]0,1[$,  by the  convexity  assumption in [H2], the exists a control $\tilde{a}$ such that
$$
\bar{\mu} ( \bH(z,w,a),\lH(z,w,a)) +(1-\bar{\mu})  ( \bH(z,w,a_1),\lH(z,w,a_1))= ( \bH(z,w,\tilde{a}),\lH(z,w,\tilde{a})) .
$$
By construction  $  \bH(z,w,\tilde{a}) \cdot \nor_1(w) =0$, therefore $\tilde{a}  \in A_0(z,w)$. Moreover,  
since $$(1-\bar{\mu}) =\frac{ \bH(z,w,a) \cdot \nor_1(w) }{ \bH(z,w,a) \cdot \nor_1(w)+\delta}$$
we have
 $$\vert 1-\bar{\mu}\vert\leq \frac 1\delta\vert \bH(z,w,a) \cdot \nor_1(w)-\bH(x,y,a) \cdot \nor_1(y)\vert\leq \frac 1\delta( \omega_{b}( x-z)+| y-w|   \max_{i=1,2}(\bar{L}_i)+M_b\vert  \nor_1(y)- \nor_1(w)\vert) $$
then
$$ 
| ( \bH(z,w,a),\lH(z,w,a)) -( \bH(z,w,\tilde{a}),\lH(z,w,\tilde{a}) | \leq   
$$
$$
\leq | (1-\bar{\mu}) || (   ( \bH(z,w,a),\lH(z,w,a)) -  ( \bH(z,w,a_1),\lH(z,w,a_1))  ) |  \leq  $$ 
$$
\leq  \frac {2\max( M_b,M_l  )}\delta( \omega_{b}( x-z)+(\max_{i=1,2}(\bar{L}_i)+M_bL_n)| y-w|) .
$$
Hence, also thanks to  \eqref{lipafix} we obtain
$$
|( \bH(x,y,a),\lH(x,y,a)) -  ( \bH(z,w,\tilde{a}),\lH(z,w,\tilde{a}))  |  \leq 
$$
$$
|( \bH(x,y,a),\lH(x,y,a)) -     ( \bH(z,w,a),\lH(z,w,a))  | + 
|( \bH(z,w,a),\lH(z,w,a))-  ( \bH(z,w,\tilde{a}),\lH(z,w,\tilde{a}))  |   \leq 
$$
$$
\leq   \frac {2\max( M_b,M_l  )}\delta( \omega_{b}( x-z)+(\max_{i=1,2}(\bar{L}_i)+M_bL_n)| y-w|)+|  \omega_{b}( x-z)+ \omega_{l}( x-z) +| y-w|   \max_{i=1,2}(\bar{L}_i,\bar{L}_{i,l}) |  =
$$
$$
= (\frac {2\max( M_b,M_l  )}\delta+1) \omega_{b}( x-z)+ \omega_{l}( x-z)+( (\frac {2\max( M_b,M_l  )}\delta)(\max_{i=1,2}(\bar{L}_i)+M_bL_n)+ \max_{i=1,2}(\bar{L}_i,\bar{L}_{i,l}))| y-w|
$$
and this concludes the proof.   \hfill $\Box$

Let us now prove some regularity properties on the tangential Hamiltonian $H_T$ (see also Lemma 7.2. in \cite {BBC2}). 
\begin{pro}
\label{propHtan} 
Assume [H0bis], [H1], 
and [H2], let $\H$ be a $W^{2,\infty}$-hypersurface. Let $(x,y), (z,w) \in \R^\N \times \H$.\\
The tangential Hamiltonian defined in \eqref{def:HamHT} satisfies the following Lipschitz properties with respect $x$ and $p_\H$:\\
There exists a constant  $M$ (which can be estimated by $M\leq M_b$) such that  for any $p_\H\in T_y\H$ and $q_\H\in T_y\H$
\begin{equation}  \label{HT-lip:pH}
\vert \HT(x,y,p_\H)-\HT(x,y,q_\H)\vert\leq M\vert p_\H-q_\H\vert .
\end{equation}
There exists a constant $C$ such that for any $p_\H\in T_y\H$ and $q_\H\in T_w\H$
\begin{equation}  \label{HT-lip:x}  
\vert \HT(x,y,p_\H)-\HT(z,w,q_\H)\vert\leq C(((\vert p_\H\vert+\vert q_\H\vert)\cdot |  \omega_{b}( x-z) +| y-w|   | + | \omega_{l}( x-z) +| y-w|   | )+\vert q_\H-p_\H\vert)
\end{equation}
\end{pro}
{\bf Proof.} The proof easily follows from Lemma   \ref{lemmacont} and standard arguments.   \hfill $\Box$

\begin{rem} \label{propHtanreg}
We remark here that the results of Lemma    \ref{lemmacont}   and Proposition  \ref{propHtan}  still hold in the case of $\HTreg$ changing the constants in  \eqref{lemmaconteq} and  \eqref{HT-lip:pH}. This can be seen as in \cite[Remark 6.7]{BBC2}. 
\end{rem}

We are finally ready to prove the stability result. 
\begin{theo}   \label{teostab}    \label{teostabp} 
Assume [H0], [H1], 
and [H2]. Fix $x,p \in \R^N$. 
Let $\rho >0$  and  $u^\rho$ ,  $v^\rho$ be respectively  sequence of  sub and supersolution of 
\be
 \rho w(y) +\HHm(x,y,Dw) =  0\mbox{ in }  \R^N .
\ee 
If $(\rho u^\rho,  u^\rho )\ds (-\mu_1, u) $ and $(\rho v^\rho,  v^\rho )\ds (-\mu_2, v)$ uniformly in $\R^N$,   then   
the function $u$ is a viscosity subsolution of  
\be
  \HHm(x,y,Du) =  \mu_1 \mbox{ in }  \R^N
\ee 
while the function $v$ is a viscosity supersolution of 
\be
   \HHm(x,y,Dv) =  \mu_2 \mbox{ in }  \R^N.
\ee
Moreover, let $\rho >0$  and  $u^\rho$  be a  sequence of  subsolution of 
\be
 \rho w(y) +\HHp(x,y,Dw) =  0\mbox{ in }  \R^N .
\ee 
If $(\rho u^\rho,  u^\rho )\ds (-\mu_1, u) $  uniformly in $\R^N$,   then   
the function $u$ is a viscosity subsolution of  
\be
  \HHp(x,y,Du) =  \mu_1 \mbox{ in }  \R^N .
\ee 
\end{theo}
{\bf Proof.}  
Since we are assuming that the convergence is uniform,  this  result can be proven following  standard arguments for  stability results 
on viscosity solutions (see, for instance \cite{Ba}). Note that the only difference with the standard result is in the  proof of the limit inequality on $\H$.  
However, the standard arguments apply  thanks to the regularity of the tangential Hamiltonian  $\H_T$   and $\HTreg$
proved in Proposition   \ref{propHtan}  and Remark  \ref{propHtanreg}, respectively. 
\hfill $\Box$
\begin{rem}
Note that in this  stability result the uniform convergence assumption can not be weakened. 
Indeed, roughly speaking, this condition is necessary to apply "separately"  the  standard argument and  pass to the limit for $H_1$  only in $\Omega_1$, 
for $H_2$ only in $\Omega_2$ and for $\HT$ only in $\H$.
\end{rem}
\begin{rem}
Note that Theorem \ref{propHeffe} and Theorem \ref{Car-Hbarra+0} hold true also if instead of [H0]
we assume [H0bis].
\end{rem}



    
    


\paragraph{Acknowledgement}
The authors  were partially funded  by the ANR  project ANR HJnet (ANR-12-BS01-0008-01)
and by the EU under the 7th Framework Programme Marie Curie
Initial Training Network ``FP7-PEOPLE-2010-ITN'', SADCO project, GA number 264735-SADCO.

\end{document}